\documentclass[a4paper,draft]{amsart}
 \usepackage[cp1251]{inputenc}
 \usepackage[russian]{babel}
 \usepackage{amssymb, amsmath}
 \theoremstyle{definition}
 \newtheorem{ddd}{Определение}[section]
 
 \theoremstyle{plain}
 \newtheorem{ttt}[ddd]{Теорема}
 \newtheorem{ccc}[ddd]{Следствие}
 \newtheorem{llll}[ddd]{Лемма}
 \theoremstyle{remark}
 \newtheorem*{rrr}{Замечание}

 \newcommand{\mdeg}{\mathrm{mdeg}}
 \newcommand{\supp}{\mathrm{supp}}

\begin{document}
   \title[Централизаторы в частично коммутативных алгебрах Ли]{Централизаторы в частично коммутативных\\ алгебрах Ли
   \protect{\footnote{Работа выполнена при поддержке РФФИ, грант
   12-01-00084}}}
   \author{Порошенко Е.\,Н.}
   \date{}
   \begin{abstract}
     Работа посвящена получению полного описания централизаторов элементов частично коммутативных алгебр
     Ли. Результат приводится в явном виде в терминах порождающих
     частично коммутативной алгебры.
   \end{abstract}
   \begin{flushright}
    УДК 512.554.33
   \end{flushright}
   \maketitle
   \section{Введение}
     Пусть
     $X$~--- конечное множество и
     $G=\langle X,E\rangle$~--- неориентированный граф без петель, множеством вершин которого
     является множество
     $X$, а множеством ребер~--- некоторое отношение на множестве
     $X$, то есть подмножество множества
     $X\times X$.  Так как граф
     $G$ неориентированный, элементами множества
     $E$ являются неупорядоченные пары, которые будут обозначаться
     $\{x,y\}$, где
     $x,y\in X$.

     \emph{Частично коммутативной алгеброй Ли} над областью целостности
     $R$ называется
     $R$-алгебра
     $\mathcal{L}_R(X;G)$ с множеством порождающих
     $X$ и множеством определяющих соотношений
     $$(x_i,x_j)=0, \text{ где } \{x_i,x_j\}\in E.$$
     (Здесь и далее
     $(g,h)$ обозначает лиевское произведение элементов
     $g$ и
     $h$).
     Граф
     $G$ называется \emph{определяющим графом} соответствующей алгебры. Для упрощения обозначений,
     мы также будем обозначать эту алгебру
     $\mathcal{L}(X;G)$, если не возникнет неоднозначности.

     Таким образом, определение частично коммутативных алгебр Ли
     аналогично определениям других частично коммутативных
     структур: групп, моноидов и т.д. (см.
\cite{DK93}). Если частично коммутативные группы
     являются объектом пристального внимания (см.,например,
\cite{GT09,Ti10,Ti11,She05,She06,DKR07,Se89}), то частично
     коммутативные алгебры (как ассоциативные, так и алгебры Ли) изучены гораздо меньше.
     Однако и для этих объектов получен ряд результатов. Например,
     в
\cite{HKMLNR80} доказано, что ассоциативные частично коммутативные
     алгебры над произвольным полем изоморфны тогда и только
     тогда, когда изоморфны определяющие их графы. В
\cite{DK92} приводится алгоритм, позволяющий находить базис любой
     частично коммутативной алгеры Ли. Однако, этот алгоритм
     основан на разбиении множества
     $X$ вершин графа
     $G$ на два подмножества, одно из которых является
     независимым, поэтому данный алгоритм существенно зависит от
     структуры графа
     $G$, а значит не дает явной структуры базисов частично
     коммутативных алгебр Ли в общем случае. Явное же описание базисов
     частично коммутативных алгебр Ли получено в
\cite{por11}.

     Целью данной работы является полное описание централизаторов
     элементов частично коммутативных алгебр Ли. Статья начинается с
     формулировки основных определений и результатов, необходимых
     для дальнейшей работы  (см. раздел~%
\ref{1}). Далее работа построена по принципу от
     частного к общему. В разделе~%
\ref{2} описываются централизаторы порождающих алгебры Ли, то есть
     элементов множества
     $X$. Раздел~%
\ref{3} посвящен нахождению централизаторов линейных комбинаций
     порождающих. Промежуточные результаты сформулированы в следствии
\ref{centlincomb} и теореме
 \ref{lincombgen}.  Наконец, в разделе~%
\ref{4} доказывается теорема~
 \ref{gencase}, дающая описание централизаторов в самом общем
     случае.

   \section{Предварительные сведения}\label{1}
     Пусть
     $\mathcal{L}_R(X;G)$~--- частично коммутативная алгебра Ли с
     определяющим графом
     $G$. Если
     $\{x,y\}\in E$, то будем писать
     $x\leftrightarrow y$. Аналогично, если
     $Y\subseteq X$ и
     $x\in X$~--- порождающий, смежный со всеми вершинами множества
     $x$ в графе
     $G$, то будем писать
     $x\leftrightarrow Y$. Наконец, такой же смысл придадим записи
     $Y \leftrightarrow Z$ для
     $Y,Z\subseteq X$: каждый порождающий множества
     $Y$ соединен в графе
     $G$ ребром с каждым порождающим множества
     $Z$. В частности, так как граф
     $G$ не содержит петель, из
     $x\leftrightarrow Y$ следует, что
     $x\not \in Y$, а
     $Y\leftrightarrow Z$ влечет
     $Y\cap Z=\varnothing$. Данное обозначение будем использовать как глобальное, то есть оно будет означать,
     что вершины (или множества вершин) соединены
     ребрами в определяющем графе частично коммутативной алгебры, а не в графе, рассматриваемом
     в тот или иной момент.

     Нам потребуются еще некоторые обозначения, связанные с графами. Пусть
     $H$~--- произвольный неориентированный граф, с множеством
     вершин. Через
     $V(H)$ и
     $E(H)$ будем обозначать соответственно множество вершин и
     множество ребер этого графа. Далее, пусть
     $V'\subseteq V(H)$. Через
     $H(V')$ обозначим подграф графа
     $H$, построенный на множестве вершин
     $V'$. Наконец, напомним, что через
     $\overline{H}$ обозначается граф
     $\langle V(H), (V(H))^2\backslash (E(H)\cup \mathrm{id}_{V(H)})\rangle$,
     где для произвольного множества
     $S$ через
     $\mathrm{id}_S$ обозначается тождественное бинарное отношение:
     $\mathrm{id}_S=\{(x,x)\,|\, x\in S\}$; этот граф называется дополнением графа
     $H$.

     Итак, пусть
     $X^*$~--- множество всех ассоциативных слов алфавита
     $X$ (включая пустое слово, которое мы будем обозначать
     $1$). Введем на множестве
     $X$ линейный порядок и распространим его до так называемого \emph{лексикографического} порядка
     на множестве
     $X^*$:
     \begin{enumerate}
       \item
         $u<1$ для любого непустого слова
         $u$.
       \item
         по индукции,
         $u<v$, если
         $u=x_iu'$,
         $v=x_jv'$, где
         $x_i,x_j\in X$ и
         $x_i<x_j$ или
         $x_i=x_j$ и
         $u'<v'$.
     \end{enumerate}

     \begin{ddd}\label{alsw}
       Слово
       $u\in X^*$ называется \emph{ассоциативным словом Линдона~--- Ширшова}, если для
       любых непустых
       $v,w\in X^*$, таких что
       $u=vw$, выполнено
       $wv<u$.
     \end{ddd}

     Также рассмотрим множество всех неассоциативных мономов на
     множестве порождающих
     $X$ (здесь мы исключаем пустое слово из рассмотрения), то есть
     множество всех слов со всеми возможными способами расставления на них скобок (обозначим его
     $X^{+}$). В следующем определении если
     $[u]$~--- произвольный неассоциативный моном, то будем обозначать через
     $u$ ассоциативное слово, полученное из
     $[u]$ стиранием скобок.

     \begin{ddd}\label{nalsw}
       \renewcommand{\theenumi}{\roman{enumi}}
       Cлово
       $[u]\in X^{+}$ называется \emph{неассоциативным словом Линдона~--- Ширшова}, если
       \begin{enumerate}
         \item
           Ассоциативное слово
           $u$ является ассоциативным
           словом Линдона~--- Ширшова;
         \item
           Если
           $[u]=([u_1],[u_2])$, то
           $[u_1]$ и
           $[u_2]$~--- неассоциативные слова Линдона~--- Ширшова и
           $u_1>u_2$.
         \item
           Если
           $[u_1]=([u_{11}],[u_{12}])$, то
           $u_2\geqslant u_{12}$.
       \end{enumerate}
     \end{ddd}

     Слова данного вида были впервые рассмотрены в
\cite{sh53}, а их аналог для групп~--- в
 \cite{CFL58}.

     Следует отметить, что множества ассоциативных и неассоциативных слов Линдона~--- Ширшова
     зависят от порядка на множестве порождающих
     $X$. Соответственно, множества ассоциативных и неассоциативных слов Линдона~--- Ширшова на множестве
     порождающих
     $X$ при его фиксированном упорядочении
     $O$ будем обозначать
     $LSA(X;O)$ и
     $LS(X;O)$, а элементы этих множеств будем называть
     $O$-LSA- и
     $O$-LS-словами. Если не возникает неоднозначности, можно
     вместо
     $LSA(X;O)$ и
     $LS(X;O)$ писать
     $LSA(X)$ и
     $LS(X)$, а
     $O$-LSA- и
     $O$-LS слова будем называть просто LSA- и
     LS-словами.

     Как следует из
\cite{sh58}, множество
     $LS(X)$ образует базис свободной алгебры Ли c
     множеством порождающих
     $X$ (будем обозначать эту алгебру
     $\mathcal{L}_R(X)$ или
     $\mathcal{L}(X)$, если из контекста ясно о какой области целостности
     $R$ идет речь).

     В работе
\cite{sh58} также было показано, что на любом
     LSA-слове можно единственным образом расставить
     скобки, чтобы получить LS-слово: для любого
     $u\in LS(X)$, такого что
     $\ell(u)\geqslant 2$ имеем
     $[u]=([v],[w])$, где
     $w$~--- наибольший собственный суффикс слова
     $u$, являющийся
     LSA-словом (соответственно, для
     $v$ и
     $w$ можно применить индукцию по длине слова). Иными словами, существует биективное отображение
     $[\,\cdot\,]: LSA(X)\rightarrow LS(X)$. В дальнейшем, для
     любого слова
     $u\in LSA(X)$ мы будем обозначать его образ под действием этого
     отображения через
     $[u]$. Наконец, для
     $[u],[v]\in LS(X)$ будем писать
     $[u]<[v]$, если
     $u<v$. Аналогичный смысл будем придавать выражениям
     $[u]\leqslant [v]$,
     $[u]>[v]$ и
     $[u]\geqslant [v]$.

     Напомним два свойства слов Линдона~--- Ширшова.
     \begin{llll}[\cite{BC08}, лемма~2.12] \label{LSAL}
       Если
       $u,v\in LSA(X)$ и
       $u>v$, то
       $uv\in LSA(X)$.
     \end{llll}

     \begin{llll}\label{less}
       Пусть
       $[u],[v]\in LS(X)$ и
       $[u]>[v]$, тогда в алгебре
       $\mathcal{L}_{R}(X)$ имеет место представление
       $([u],[v])=[uv]+\sum_{i}\alpha_i [w_i]$, где
       $\alpha_i\in R\backslash\{0\}$, все
       $[w_i]\in LS(X)$ и
       $w_i<uv$.
     \end{llll}

     Лемма~%
\ref{less} с очевидностью
   следует из
\cite{sh58} (леммы~2 и 3) и леммы
 \ref{LSAL}.

     Для частично коммутативной алгебры Ли
     $\mathcal{L}(X;G)$ с множеством порождающих
     $X$ определим по индукции \emph{частично коммутативные слова
     Линдона~--- Ширшова}:
     \begin{enumerate}
       \item
         Элементы множества
         $X$ являются PCLS-словами.
       \item
         LS-слово
         $[u]$, такое что
         $\ell([u])>1$, является PCLS-словом, если
         $[u]=([v],[w])$, где
         $[v]$ и
         $[w]$~--- PCLS-слова, причем в графе
         $G$ первая буква лиевского монома
         $[w]$ не соединена ребром хотя бы с одной буквой монома
         $[v]$.
     \end{enumerate}
     Так как граф
     $G$ без петель, в частности, PCLS-словами являются мономы вида
     $[u]=([v],[w])$, в которых первая буква
     $[w]$ принадлежит
     $\supp([v])$.

     В
\cite{por11} было доказано, что множество  частично
     коммутативных слов Линдона~--- Ширшова образует
     базис алгебры
     $\mathcal{L}(X;G)$. Это множество, как и множества
     LSA- и LS-слов, определяется упорядочением множества порождающих
     $X$, соответственно, мы будем использовать аналогичные
     обозначения: множество всех частично коммутативных слов при
     упорядочении
     $O$ множества порождающих
     $X$ будем обозначать
     $PCLS(X;O)$, а элементы этого множества будем называть
     $O$-PCLS-словами. Если упорядочение
     $O$ определено контекстом, то, аналогично ранее введенным обозначениям, будем опускать в
     явное указание на упорядочение порождающих.

     Элементы множеств
     $LS(X;O)$ и
     $PCLS(X;O)$ будем также называть
     $O$-\emph{базисными мономами} алгебр
     $\mathcal{L}(X)$ и
     $\mathcal{L}(X;G)$ соответственно или \emph{базисными мономами относительно упорядочения}
     $O$ в алгебре
     $\mathcal{L}(X)$ (соответственно
     $\mathcal{L}(X;G)$). Опять же, если упорядочение
     $O$ множества
     $X$ определено контекстом, можно убирать явное указание на упорядочение  и называть
     соответствующие мономы базисными. Если это необходимо, будем явно
     указывать, что используется упорядочение, порожденное порядком
     $O$ на множестве
     $X$: вместо
     ``$<$'' будем писать
     ``$<_O$'' и так далее.

     С этого места будем считать, что
     $X=\{x_1,x_2, \dots, x_n\}$.
      \begin{ddd}
       Пусть
       $u$~--- лиевский моном на множестве порождающих
       $X$.\emph{Мультистепенью} монома
       $u$ назовем вектор
       $\overline{\delta}=(\delta_{1},\delta_2,\dots,\delta_n)$, где
       $\delta_i$~--- это число вхождений порождающего
       $x_i$ в моном
       $[u]$.
     \end{ddd}
     Мультистепень монома
     $[u]$ будем обозначать через
     $\mdeg([u])$, а число вхождений буквы
     $x_i$ в моном
     $[u]$~--- через
     $\mdeg_i ([u])$.

     Пусть
     $[u]$~--- лиевский моном и
     $(\delta_1,\dots,\delta_n)$~--- его мультистепень. Через
     $\supp([u])$ будем обозначать множество
     $\{x_i\,|\, \delta_i \neq 0 \}$.

     \begin{ddd}
       Ненулевой элемент
       $g$ алгебры
       $\mathcal{L}(X;G)$ называется \emph{однородным}, если
       $g$ представим в виде линейной комбинации лиевских мономов одной и той же мультистепени
       $\overline{\delta}=(\delta_{1},\delta_2,\dots,\delta_n)$.
     \end{ddd}

     Отметим, что в силу однородности тождеств алгебры Ли и определяющих соотношений частично коммутативной
     алгебры справедливо следующее утверждение. Если
     $0=\sum_{i}g_i$, где каждое слагаемое
     $g_i$ является однородным лиевским многочленом, причем все
     $g_i$ имеют различные мультистепени, то
     $g_i=0$ для всех
     $i$. В частности, если элемент
     $g$ однородный, то при любом упорядочении
     $O$ множества
     $X$  все мономы разложения
     $g$ в линейную комбинацию
     $O$-базисных мономов имеют одну и ту же
     мультистепень. Более того, мультистепень мономов, входящих в разложение
     однородных элементов не зависит от
     конкретного упорядочения множества порождающих.
     Соответственно, можно распространить понятие мультистепени на множество однородных элементов алгебры
     $\mathcal{L}(X;G)$, полагая
     $\mdeg(g)=\mdeg([u])$, где
     $[u]$~--- произвольный базисный лиевский моном, входящий в разложение
     $g$ (при произвольном упорядочении порождающих).

     Так как мультистепень является вектором, при работе с мультистепенями можно
     использовать те же обозначения, что и при работе с векторами. В частности, удобно использовать понятие суммы
     мультистепеней, определив ее как сумму векторов, соответствующих исходным
     мультистепеням. Легко видеть, что в силу однородности тождеств алгебры Ли и
     определяющих соотношений частично коммутативной алгебры Ли,
     если
     $(g,h)\neq 0$ для однородных элементов
     $g,h\in \mathcal{L}(X;G)$, то
     $(g,h)$ также является однородным элементом и
     $\mdeg((g,h))=\mdeg(g)+\mdeg(h)$.
     Для произвольной мультистепени
     $\overline{\delta}=(\delta_1,\delta_2,\dots, \delta_n)$ ее
     $i$-ую координату будем называть \emph{степенью
     $\overline{\delta}$ относительно
     $x_i$}.

     Для произвольной перестановки
     $\sigma\in S_n$ определим порядок на множестве мультистепеней: будем говорить, что
     мультистепень
     $\overline{\delta}=(\delta_1,\delta_2,\dots,\delta_n)$ \emph{больше} мультистепени
     $\overline{\xi}=(\xi_1,\xi_2\dots,\xi_n)$ и писать
     $\overline{\delta}>\overline{\xi}$, если для некоторого
     $k\in\{1,2,\dots,n\}$ выполнено
     $\delta_{\sigma(i)}=\xi_{\sigma(i)}$ для
     $i>k$ и
     $\delta_{\sigma(k)}>\xi_{\sigma(k)}$. Этот порядок будем
     обозначать
     $\Omega(\sigma)$ или
     $\Omega(x_{\sigma(n)},\dots, x_{\sigma(1)})$. Иными словами, мы вначале
     сравниваем степени
     $\overline{\delta}$ относительно
     $x_{\sigma(n)}$, затем, в случае их равенства, степени относительно
     $x_{\sigma(n-1)}$ и так далее. Посредством выбора подстановки
     $\sigma$ можно задать любой порядок, в котором будут
     сравниваться координаты векторов, задающих мультистепени.

     Пусть
     $y_1,\dots y_m \in X$. Через
     $\Omega(y_m,y_{m-1}\dots, y_1)$ будем множество порядков, в которых
      у мультистепеней вначале сравниваются степени относительно порождающего
     $y_m$, затем степени относительно породждающего
     $y_{m-1}$ и так далее. Это обозначение мы будем
     использовать, когда нам важен лишь порядок сравнения
     степеней относительно порождающих из некоторого подмножества
     $Y\subseteq X$. Оно хорошо согласуется с обозначениями конкретного порядка
     $\Omega(\sigma)$ и
     $\Omega(x_{\sigma(n)},\dots, x_{\sigma(1)})$ если рассматривать последние также
     как множество порядков, которое в этом случае будет одноэлементным.
     Если нам необходимо подчеркнуть, что мультистепени сравниваются смысле
     конкретного упорядочения
     $\mathcal{O}$, то будем использовать обозначения
     $<_{\mathcal{O}}$,
     $\leqslant_{\mathcal{O}}$,
     $>_{\mathcal{O}}$ и
     $\geqslant_{\mathcal{O}}$ вместо
     $<$,
     $\leqslant$,
     $>$ и
     $\geqslant$ соответственно.

     Отметим, что в конкретном случае
     порядок мультистепенях не обязательно связан с порядком на
     порождающих, то есть не обязательно вначале сравниваются
     степени относительно наибольшего порождающего, затем относительно второго по
     порядку, и так далее.

     Очевидно, что
     $\supp([u])$ принимает одно и то же значение для всех базисных мономов
     $[u]$, входящих в разложение однородного элемента
     $g$ в виде линейной комбинации базисных элементов. Следовательно, можно по определению положить
     $\supp (g)=\supp([u])$, где
     $[u]$~--- произвольный
     $O$-базисный моном, входящий в разложение элемента
     $g$ (с ненулевым коэффициентом) при любом упорядочении
     $O$. Пусть элемент
     $g$ не является однородным, тогда имеет место представление
     $g=\sum_{\overline{\delta}}g_{\overline{\delta}}$ в виде
     суммы однородных слагаемых различных мультистепеней. В этом
     случае положим
     $\supp(g)=\bigcup\limits_{\overline{\delta}}\supp(g_{\overline{\delta}})$.

     \begin{ddd}
       Пусть
       $g$~--- произвольный элемент алгебры частично коммутативной
       алгебры Ли
       $\mathcal{L}(X;G)$ и пусть
       $C(g)$~--- централизатор элемента
       $g$. \emph{Централизатором  в коммутанте} называется
       множество
       $$\mathcal{C}(g)=C(g)\cap \mathcal{L}'(X;G),$$
       где
       $\mathcal{L}'(X;G)$ это производная алгебры
       $\mathcal{L}(X;G)$.
     \end{ddd}

     \begin{ddd}
       Пусть
       $\mathfrak{L}$~--- алгебра Ли над областью целостности
       $R$. Элементы
       $g,h\in \mathfrak{L}$ называются \emph{пропорциональными},
       если для некоторых
       $\lambda,\mu \in R$ выполнено
       $\lambda g=\mu h$. Если
       $g$ и
       $h$ пропорциональны, то будем писать
       $g\backsim h$.
     \end{ddd}

   \section{Централизаторы порождающих частично коммутативных алгебр}\label{2}
     Пусть
     $x$~--- некоторый порождающий. Для нахождения
     $C(x)$ нам потребуются две вспомогательные леммы.
     \begin{llll}\label{xinsuppg}
       Пусть
       $g$ ненулевой однородный элемент алгебры
       $\mathcal{L}_{R}(X;G)$, не являющийся порождающим. Если
       $x\in \supp(g)$, то
       $[g,x]\neq 0$.
     \end{llll}
     \begin{proof}
       Пусть
       $g$~--- ненулевой элемент алгебры
       $\mathcal{L}_R(X;G)$ и пусть
       $O$~--- порядок на множестве порождающих
       $X$, при котором
       $x$~--- наибольший порождающий. При разложении
       $g$ по базису, заданному порядком
       $O$ получаем
       $g=\sum_{i}\alpha_i [u_i]$, где
       $[u_i]$~--- различные
       $O$-PCLS-слова, причем все мономы
       $[u_i]$ имеют один и тот же состав и
       $x\in \supp(g)$, а
       $\alpha_i$~--- ненулевые элементы области целостности
       $R$.

       Так как
       $x$~--- наибольший порождающий, эта буква является первой в записи всех
       $[u_i]$. Имеем
       $([u_i],x)=-(x,[u_i])$. Причем
       $(x,[u_i])$~--- PCLS-слово. Действительно, так как
       $x$~--- наибольший порождающий, моном
       $(x,[u_i])$ является словом Линдона~--- Ширшова. Кроме того,
       $[u_i]$ является PCLS-словом, очевидно, что и
       $x$~--- тоже PCLS-слово. Значит, так как первая буква
       $[u_i]$ также равна
       $x$, получаем
       $(x,[u_i])$ тоже PCLS-слово. Следовательно, имеем
       $$(g,x)=\sum_i (-\alpha_i)(x,[u_i]).$$
       Так как все мономы
       $(x,[u_i])$ являются
       $O$-базисными и различны, а все
       $\alpha_i$ не равны
       $0$, получаем
       $(g,x)\neq 0$.
     \end{proof}
     \begin{llll}\label{xnotinsuppg}
       Пусть
       $g$ ненулевой однородный элемент алгебры
       $\mathcal{L}_R(X;G)$. Если вершина
       $x\not \in \supp(g)$ и
       $x \nleftrightarrow\supp(g)$, то
       $(g,x)\neq 0$.
     \end{llll}
     \begin{proof}
       Пусть
       $g$~--- ненулевой элемент алгебры
       $\mathcal{L}_R(X;G)$, такой что
       $x\not \in \supp(g)$ и пусть
       $y\in \supp(g)$ вершина, не смежная с вершиной
       $x$ в графе
       $G$. Рассмотрим порядок
       $O$, в котором
       $x$~--- наибольший  порождающий, а
       $y$~--- второй по старшинству. При разложении
       $g$ по базису, заданному порядком
       $O$, получаем
       $g=\sum_{i}\alpha_i [u_i]$, где
       $[u_i]$~--- различные PCLS-слова относительно порядка
       $O$, причем все мономы
       $[u_i]$ имеют одну и ту же мультистепень и начинаются с
       порождающего
       $y$, а
       $\alpha_i$~--- ненулевые элементы области целостности
       $R$.

       Имеем
       $([u_i],x)=-(x,[u_i])$. Причем
       $(x,[u_i])$~---PCLS-слово. Действительно, так как
       $x$~--- наибольший порождающий, моном
       $(x,[u_i])$ является словом Линдона~--- Ширшова. Кроме того,
       $[u_i]$ является PCLS-словом, очевидно, что и
       $x$~--- тоже PCLS-слово. Значит, так как первая буква
       $[u_i]$, не смежна с
       $x$, получаем
       $(x,[u_i])$ тоже PCLS-слово. Следовательно, имеем
       $$(g,x)=\sum_i (-\alpha_i)(x,[u_i]).$$
       Так как все мономы
       $(x,[u_i])$ являются
       $O$-базисными и различны, а все
       $\alpha_i$ не равны
       $0$, получаем
       $(g,x)\neq 0$.
     \end{proof}

     \begin{ttt}\label{centgen}
       В частично коммутативной алгебре Ли
       $\mathcal{L}_R(X;G)$ централизатор в коммутанте для
       порождающго
       $x$ состоит из всех элементов алгебры, представимых в виде
       $\sum_i \alpha_i [u_i]$, где
       $[u_i]$~--- базисные лиевские мономы, такие что
       $\ell(u_i)\geqslant 2$ и
       $x \leftrightarrow \supp([u_i])$.
     \end{ttt}
     \begin{proof}
       Очевидно, что все элементы, указанного вида лежат в
       централизаторе. Действительно, пусть
       $[u_i]=(y_1,y_2,\dots y_k)$, с произвольной расстановкой
       скобок. Тогда
       $([u_i],x)=((y_1,y_2,\dots y_k),x)=\sum_{i=1}^k (y_1,\dots, (y_i,x),\dots, y_k)$.
       Так как каждое слагаемое содержит множитель
       $(y_i,x)=0$, получаем
       $([u_i],x)=0$.

       Докажем, что других элементов в централизаторе нет.
       Действительно, пусть
       $g$~--- элемент централизатора. Имеет место представление
       $g=\sum_i g_i$, где
       $g_i$~--- однородные лиевские многочлены различной
       мультистепени. Тогда
       $0=(g,x)=\sum_i (g_i,x)$, причем
       $(g_i,x)$~--- также однородные элементы и те из них, которые не равны нулю, имеют различные
       мультистепени. Отсюда получаем, что
       $(g_i,x)=0$, то есть
       $g_i$~--- также элементы из централизатора и из лемм~%
\ref{xinsuppg} и
 \ref{xnotinsuppg} следует, что
       $x\not \in \supp(g_i)$ и
       $x\leftrightarrow \supp(g_i)$, что и требовалось.
     \end{proof}
     Из теоремы
\ref{centgen} получаем следующее утверждение.
     \begin{ccc}\label{centlincomb}
       В частично коммутативной алгебре Ли
       $\mathcal{L}_R(X;G)$ централизатор для порождающего
       $x$ состоит из всех элементов алгебры, представимых в виде
       $\alpha x+\sum_i \alpha_i [u_i]$, где
       $u_i$~--- базисные лиевские мономы, такие что
       $x\leftrightarrow \supp([u_i])$.
     \end{ccc}

   \section{Централизаторы линейных комбинаций порождающих}
   \label{3}

     Мы начнем этот раздел с нахождения централизатора в коммутанте для произвольной линейной комбинации порождающих
     $\sum_{j=1}^m \alpha_{i_j} x_{i_j}$, где все
     $\alpha_{i_j}\neq 0$.

     \begin{ttt}\label{centinter}
       Пусть
       $\mathcal{L}(X;G)$~--- частично коммутативная алгебра Ли над областью целостности
       $R$, где
       $X=\{x_1,x_2,\dots,x_n\}$. Тогда для любого набора различных порождающих
       $x_{i_1},x_{i_2},\dots, x_{i_m}$ и для любых
       $\alpha_{i_1},\alpha_{i_2},\dots,\alpha_{i_m}\in R\backslash\{0\}$
       выполнено
       $$\mathcal{C}\bigl(\sum_{j=1}^m \alpha_{i_j} x_{i_j}\bigr)=\bigcap_{j=1}^{m} \mathcal{C}(x_{i_j}).$$
     \end{ttt}
     \begin{proof}
       Пусть
       $y_1,y_2, \dots, y_k$ некоторые различные порождающие алгебры
       $\mathcal{L}(X;G)$. Очевидно, что для любых
       $\alpha_1,\alpha_2, \dots, \alpha_k \in R\backslash\{0\}$
       имеет место включение
       $\bigcap_{i=1}^k \mathcal{C}(y_i)\subseteq \mathcal{C}(\sum_{i=1}^k \alpha_i y_i)$.

       Докажем обратное включение. Пусть
       $O$ порядок на порождающих, в котором
       $y_k>y_{k-1}>\dots >y_1$, а остальные порождающие меньше, чем
       $y_1$ и пусть
       $\mathcal{O}\in \Omega(y_k,y_{k-1},\dots, y_1)$.
       Обозначим
       $\sum_{i=1}^k \alpha_i y_i$ через
       $h$.
       Пусть
       $g\in \mathcal{C}(h)$. Имеет место
       разложение
       \begin{equation} \label{decomp}
         g=\sum_{\overline{\delta}} g_{\overline{\delta}},
       \end{equation}
       где
       $g_{\overline{\delta}}$~--- однородная компонента
       мультистепени
       $\overline{\delta}$ лиевского
       многочлена
       $g$.

       Для доказательства применим индукцию по
       $k$, то есть по количеству слагаемых в разложении
       $h$ в линейную комбинацию порождающих.

       Сначала докажем, что
       $g_{\overline{\delta}}\in \mathcal{C}(y_k)$ для всех
       слагаемых
       $g_{\overline{\delta}}$ из правой части равенства
(\ref{decomp}).

       Для удобства переупорядочим координаты в мультистепенях
       таким образом, чтобы последняя координата соответствовала
       порождающему
       $y_k$, предпоследняя
       $y_{k-1}$ и так далее, то есть порождающему
       $y_i$ соответствует координата под номером
       $n-k+i$.
       Пусть
       $\overline{\varepsilon}=(\varepsilon_1,\dots,\varepsilon_n)$~--- наибольшая из мультистепеней всех мономов,
       входящих в разложение
       $g$ (относительно порядка
       $\mathcal{O}$). Имеем
       $$\left(g_{\overline{\varepsilon}},\sum_{i=1}^k \alpha_i y_i\right)=\sum_{i=1}^k \alpha_i (g_{\overline{\varepsilon}},y_i).$$

       Легко видеть, что мультистепени всех базисных мономов, входящих в
       $(g,\sum_{i=1}^k \alpha_i y_i)$, но не входящих в
       $(g_{\overline{\varepsilon}},y_k)$ меньше, чем
       $(\varepsilon_1,\dots,\varepsilon_{n-1},\varepsilon_n+1)$, в то время как если
       $(g_{\overline{\varepsilon}},y_k)\neq 0$, то
       $$\mdeg((g_{\overline{\varepsilon}},y_k))=(\varepsilon_1,\dots,\varepsilon_{n-1},\varepsilon_n+1).$$
       Действительно, пусть
       $\mdeg([u])=(\delta_1,\dots,\delta_{n-1},\delta_{n})$ для
       некоторого монома
       $[u]$, входящего в разложение
       $g$. Тогда если
       $([u],y_i)\neq 0$, то
       $\mdeg(([u],y_i))=(\delta_1,\dots,\delta_{n-k+i-1},\delta_{n-k+i}+1,\dots, \delta_{n})$.
       Так как
       $\overline{\delta}\leqslant \overline{\varepsilon}$, при
       $i\neq k$ имеем
       $\delta_{n}\leqslant \varepsilon_n$, т.е.
       $\delta_{n}<\varepsilon_{n}+1$. Следовательно,
       $\mdeg(([u],y_i))<_{\mathcal{O}}(\varepsilon_1,\dots, \varepsilon_{n-1},\varepsilon_n+1)$, если
       $([u],y_i)\neq 0$. При
       $i=k$ получаем если
       $([u],y_k)\neq 0$, то
       $\mdeg((u,y_k))=(\delta_1,\dots,\delta_{n-1},\delta_{n}+1)$.
       Легко видеть, что из
       $\mdeg([u])<_{\mathcal{O}}\mdeg(g_{\overline{\varepsilon}})$ следует
       $\mdeg(([u],y_k))<_{\mathcal{O}}(\varepsilon_1,\dots, \varepsilon_{n-1},\varepsilon_n+1)$.
       Таким образом, мы получили, что если
       $(g_{\varepsilon},y_k)\neq 0$, то
       $\mdeg(([u],y_i))<_{\mathcal{O}}\mdeg((g_{\overline{\varepsilon}},y_k))$
       при
       $i \neq k$ или
       $\mdeg([u])\neq \overline{\varepsilon}$, то есть для всех мономов, входящих в разложение
       $(g,h)$, но не входящих в разложение
       $(g_{\overline{\varepsilon}},y_k)$.

       Таким образом, мономы в многочлене
       $(g_{\overline{\varepsilon}},y_k)$ не могут сократиться с
       мономами других мультистепеней. Следовательно,
       $(g_{\overline{\varepsilon}},y_k)=0$, значит, в силу леммы~%
\ref{xinsuppg},
       $y_k\not \in \supp(g_{\overline{\varepsilon}})$, то есть
       $\overline{\varepsilon}=\mdeg(g_{\overline{\varepsilon}})=(\varepsilon_1,\dots,\varepsilon_{n-1},0)$.
       Более того, для любого монома
       $[u]$ из разложения
       $g-g_{\varepsilon}$ по базису имеем
       $\mdeg_n([u])=0$, так как выполняется неравенство
       $\mdeg([u])<_{\mathcal{O}}\overline{\varepsilon}$.

       Пусть
       $\overline{\zeta}=(\zeta_1,\dots,\zeta_{n-1},0)$~--- старшая из мультистепеней базисных
       мономов, входящих в разложение
       $g-g_{\overline{\varepsilon}}$ (относительно порядка
       $\mathcal{O}$). Чтобы доказать, что
       $(g_{\overline{\zeta}},y_k)=0$, воспользуемся рассуждениями, аналогичными использованным
       выше. Действительно, если
       $(g_{\overline{\zeta}},y_k)\neq 0$, то
       $\mdeg((g_{\overline{\zeta}},y_k))=(\zeta_1,\dots,\zeta_{n-1},1)$.

       Если
       $(g_{\overline{\delta}},y_i)\neq 0$ для некоторой мультистепени
       $\overline{\delta}=(\delta_1,\dots,\delta_{n-1},0)$ и для некоторого
       $i\neq k$, то
       $$\mdeg((g_{\overline{\delta}},y_i))=(\delta_1,\dots, \delta_{n-k+i}+1,\dots, \delta_{n-1},0)<_{\mathcal{O}}(\zeta_1,\dots,\zeta_{n-1},1)$$
       для любого
       $\overline{\delta}$. Значит базисные мономы,
       входящие в запись
       $(g_{\overline{\zeta}},y_k)\neq 0$ не могут сократиться с
       базисными мономами, входящими в запись многочленов
       $(g_{\overline{\delta}},y_i)$, если
       $i\neq k$.

       Далее, если
       $(g_{\overline{\delta}},y_k)\neq 0$, то
       $$\mdeg((g_{\overline{\delta}},y_k))=(\delta_1,\dots, \delta_{n-1},1)<_{\mathcal{O}}
          (\zeta_1,\dots, \zeta_{n-1},1)=\mdeg((g_{\overline{\zeta}},y_k))$$
       для
       $\overline{\delta}<_{\mathcal{O}}\overline{\zeta}$.
       Таким образом, базисные мономы, входящие в запись
       $(g_{\overline{\zeta}},y_k)\neq 0$ не могут сократиться и с
       базисными мономами, входящими в запись многочленов
       $(g_{\overline{\delta}},y_k)$ при
       $\overline{\delta}<_{\mathcal{O}}\overline{\zeta}$. Из доказанного следует
       $(g_{\overline{\zeta}},y_k)=0$. Проводя аналогичные
       рассуждения с остальными мономами, получаем
       $(g_{\overline{\delta}},y_k)=0$ для всех
       $g_{\overline{\delta}}$ из правой части равенства
(\ref{decomp}). Таким образом, мы доказали, что
       $g_{\overline{\delta}}\in \mathcal{C}(y_k)$, следовательно
       \begin{equation}\label{simpcent}
         g \in \mathcal{C}(y_k).
       \end{equation}   Имеем
       \begin{equation} \label{simpl}
         \left(g,\sum_{i=1}^k \alpha_i y_i\right)=
         \left(g,\sum_{i=1}^{k-1} \alpha_i y_i\right)+\alpha_k(g,y_k)=
         \left(g,\sum_{i=1}^{k-1} \alpha_i y_i\right).
       \end{equation}
       Таким образом, если
       $g\in \mathcal{C}(\sum_{i=1}^k \alpha_i y_i)$, то из
(\ref{simpl}) следует, что
       $(g,\sum_{i=1}^{k-1} \alpha_i y_i)=0$, а значит
       $g\in \mathcal{C}(\sum_{i=1}^{k-1} \alpha_i y_i)$. По
       индукционному предположению
       $g\in \bigcap_{i=1}^{k-1} \mathcal{C}(y_i)$. Отсюда и из
(\ref{simpcent}) получаем
       $g\in \bigcap_{i=1}^{k} \mathcal{C}(y_i)$, что и требовалось
       доказать.
     \end{proof}
     \begin{rrr}
       Из только что доказанной теоремы следует, что
       $\mathcal{C}\left(\sum_{i=1}^k \alpha_i y_i\right)$ состоит
       из всех многочленов
       $g$, представимых в виде
       $g=\sum_{j} \beta_j [u_j]$, где
       $\ell([u_j])\geqslant 2$ для всех
       $[u_j]$, а также для каждого монома
       $[u_j]$ и для всех
       $i$ выполнено
       $y_i \leftrightarrow \supp([u_j])$.
     \end{rrr}

     Осталось описать линейные комбинации порождающих, лежащие в
     централизаторе данного элемента (также линейной комбинации
     порождающих).

     \begin{ddd}\label{proporc}
       Пусть
       $R$~--- область целостности. Упорядоченные наборы из
       $k$ элементов
       $(\alpha_1,\dots,\alpha_k),(\beta_1,\dots, \beta_k)$
       называются \emph{пропорциональными}, если существуют элементы
       $\lambda,\mu \in R$, хотя бы один из которых отличен от нуля, такие что
       $\lambda\alpha_i=\mu\beta_i$ для каждого
       $i=1,2\dots,k$.
     \end{ddd}

     Очевидно, что если
     $g=\sum_{i=1}^k \alpha_i [u_i]$ и
     $h=\sum_{i=1}^k \beta_j [u_i]$~--- разложения лиевских
     многочленов
     $g$ и
     $h$ в линейные комбинации базисных мономов (по произвольному базису) в
     произвольной алгебре Ли
     $\mathfrak{L}$, то
     эти многочлены пропорциональны тогда и только тогда, когда
     пропорциональны наборы
     $(\alpha_1,\dots, \alpha_k)$ и
     $(\beta_1,\dots, \beta_k)$.

     \begin{llll} \label{proporl}
       Пусть
       $A=(\alpha_1,\dots,\alpha_k)$ и
       $B=(\beta_1,\dots, \beta_k)$~--- упорядоченные наборы элементов из области целостности
       $R$, причем
       $\alpha_i \neq 0$ для всех
       $i$. Если
       $A$ и
       $B$ пропорциональны, то
       $\alpha_i\beta_j=\alpha_j\beta_i$ для любых
       $i,j=1, 2 \dots k$. Обратно, если для некоторого
       $s$ выполнено
       $\alpha_s\beta_i=\alpha_i\beta_s$ для любого
       $i=1, 2, \dots, k$, то наборы
       $A$ и
       $B$ пропорциональны.
     \end{llll}
     \begin{proof}
       Пусть
       $A$ и
       $B$ пропорциональны. Если
       $\lambda=0$, то, так как
       $\mu\neq 0$, имеем
       $\beta_1=\beta_2=\dots=\beta_k=0$, следовательно,
       $\alpha_i\beta_j=0$ для любых
       $i,j=1, 2 \dots k$, а значит
       $\alpha_i\beta_j=\beta_i\alpha_j$.

       Если
       $\lambda\neq 0$, то для любых
       $i$ и
       $j$ имеют место равенства
       $\lambda\alpha_i=\mu\beta_i\neq 0$,
       $\mu\beta_j=\lambda\alpha_j\neq 0$, откуда, в частности, следует, что
       $\mu\neq 0$. Перемножая эти равенства, получаем
       $(\lambda\alpha_i)(\mu\beta_j)=(\mu\beta_i)(\lambda\alpha_j)$. Таким образом,
       $0=(\lambda\alpha_i)(\mu\beta_j)-(\mu\beta_i)(\lambda\alpha_j)=
         \lambda\mu\alpha_i\beta_j-\lambda\mu\alpha_j\beta_i=\lambda\mu(\alpha_i\beta_j-\alpha_j\beta_i)$.
       Так как
       $R$~--- область целостности, и
       $\lambda,\mu\neq 0$, получаем
       $\alpha_i\beta_j-\alpha_j\beta_i=0$, что и требовалось.

       Пусть теперь для некоторого
       $s$ справедливы равенства
       $\alpha_s\beta_i=\alpha_i\beta_s$ для всех
       $i\in \{1, 2, \dots, k\}$. Если
       $\beta_s=0$, то и все остальные
       $\beta_i$ равны
       $0$, так как
       $\alpha_s\beta_i=0$,
       $\alpha_s\neq 0$ и
       $R$~--- область целостности. Полагая
       $\lambda$ равным
       $0$, а
       $\mu$ равным любому элементу из
       $R\backslash\{0\}$, получаем
       $\lambda\alpha_i=0=\mu\beta_i$,  для любого
       $i=1,2,\dots, k$, то есть
       $A$ и
       $B$~--- пропорциональны.

       Если
       $\beta_s\neq 0$ то из
       $\alpha_s\beta_i=\beta_s\alpha_i$ следует, что
       $A$ и
       $B$ удовлетворяют условию определения~%
\ref{proporc} для
       $\lambda=\beta_s$ и
       $\mu=\alpha_s$, а значит
       $A$ и
       $B$ пропорциональны.
     \end{proof}
     \begin{llll}\label{proporclemma}
       Пусть
       $\mathcal{L}(X;G)$~--- частично коммутативная алгебра Ли над
       областью целостности
       $R$ с определяющим графом
       $G$ и пусть
       $g=\sum_{i=1}^k \alpha_i y_i$ и
       $h=\sum_{j=1}^k \beta_i y_i$~--- линейные комбинации порождающих алгебры
       $\mathcal{L}(X;G)$, такие что
       $\alpha_i\neq 0$ для всех
       $i=1,2,\dots, k$, и пусть
       $(g,h)=0$. Пусть
       $H$~--- подграф графа
       $G$, порожденный множеством вершин
       $\{y_1,\dots, y_k\}$
       Тогда если в графе
       $\overline{H}$ есть путь, соединяющий вершины
       $y_s$ и
       $y_t$, то
       $\alpha_s\beta_t=\alpha_t\beta_s$.
     \end{llll}
     \begin{proof}
       Применим индукцию по длине пути, соединяющего
       $y_s$ и
       $y_t$  в графе
       $\overline{H}$.

       Пусть вершины
       $y_s$ и
       $y_t$ смежны в графе
       $\overline{H}$ (то есть не являются смежными в графе
       $H$ и, следовательно, в графе
       $G$). Элементы
       $g$ и
       $h$ могут быть записаны в виде
       $g=\alpha_s x_s+\alpha_t x_t+g'$ и
       $h=\beta_s x_s+\beta_t x_t+h'$, где
       $g'$ и
       $h'$~--- линейные комбинации порождающих из множества
       $\{y_i\,|\, i\neq s, i\neq t\}$. Имеем
       \begin{equation*}
         \begin{aligned}
           0= (g,h)= & (\alpha_s y_s+\alpha_t y_t+g',\beta_s y_s+\beta_t y_t+h') =\\
                = & (\alpha_s\beta_t-\beta_s\alpha_t) (y_s,y_t)+ (\alpha_s y_s+\alpha_t y_t,h')+
                (g',\beta_s y_s+\beta_t y_t)+(g',h')
         \end{aligned}
       \end{equation*}
       Из однородности соотношений и тождеств в частично
       коммутативных алгебрах Ли получаем
       $(\alpha_s\beta_t-\alpha_t\beta_s) (y_s,y_t)=0$. По сделанному
       предположению,
       $\{y_s,y_t\}\not \in G$, следовательно
       $\alpha_s\beta_t - \alpha_t\beta_s=0$, что и требовалось.

       Пусть утверждение леммы справедливо для всех вершин, соединенных
       путем длины меньшей, чем
       $m$, в графе
       $\overline{H}$ и пусть
       $y_s=y_{i_0},y_{i_1},\dots,y_{i_m}=y_t$~--- путь длины
       $k$, соединяющий вершины
       $y_s$ и
       $y_t$ в графе
       $\overline{H}$. Так как вершины
       $y_{i_{m-1}}$ и
       $y_t$ смежны, имеем
       $\alpha_{i_{m-1}}\beta_t=\alpha_{t}\beta_{i_{m-1}}$. Кроме того, по
       предположению индукции
       $\alpha_{i_{m-1}}\beta_{s}=\alpha_{s}\beta_{i_{m-1}}$. В силу
       леммы
\ref{proporl} упорядоченные наборы
       $(\alpha_{s},\alpha_{i_{m-1}},\alpha_{t})$ и
       $(\beta_{s},\beta_{i_{m-1}},\beta_{t})$ пропорциональны.
       Следовательно, по той же лемме
       $\alpha_s\beta_t=\alpha_t\beta_s$, что и требовалось доказать.
     \end{proof}

     Пусть
     $g=\sum_{i=1}^k \alpha_i y_i$, где
     $\alpha_i\neq 0$ для
     $i=1,2,\dots, k$ и пусть
     $h$~--- элемент, являющиеся линейной комбинацией порождающих алгебры
     $\mathcal{L}(X;G)$, причем
     $(g,h)=0$. Имеет место представление
     $h=h'+h''$, где
     $h'=\sum_{s=1}^{n-k}\gamma_s z_s$, где
    $z_s$~--- порождающие, не принадлежащие множеству
     $\{y_1,\dots y_k\}$, а
     $h''=\sum_{j=1}^k \beta_j y_j$. Из однородности отношений и тождеств частично коммутативных алгебр Ли следует, что
     $(g,h')=(g,h'')=0$.

     Так как между двумя вершинами  есть путь тогда и только тогда, когда
     эти вершины лежат в одной компоненте связности соответствующего графа, из леммы
\ref{proporclemma} следует, что если
     $(g,h'')=0$, то в
     $g$ и
     $h''$ коэффициенты при
     $y_i$-ых, лежащих в одной компоненте связности графа
     $\overline{H}$, образуют пропорциональные упорядоченные наборы.

     Таким образом, если
     $H_1,\dots H_r$~--- компоненты связности графа
     $\overline{H}$, то
     $g=\sum_{p=1}^r g_p$, где
     $g_p$~--- линейные комбинация порождающих из компоненты
     связности
     $H_p$ (с ненулевыми коэффициентами), a
     $h''=\sum_{p=1}^r h_p$, где
     $h_p$~--- линейная комбинация порождающих из компоненты
     $H_p$, при этом для каждого
     $p=1,2,\dots, r$ выполнено
     $g_p\backsim h_p$.

     Далее, получаем
     \begin{equation*}
         0= (g, h')=
           \left(\sum_{i=1}^k \alpha_i y_i, \sum_{s=1}^{n-k}\gamma_s z_s\right)=
           \sum_{i=1}^k\sum_{s=1}^{n-k}\alpha_i\gamma_s (y_i,z_s).
     \end{equation*}
     Следовательно,
     $\alpha_i\gamma_s (y_i,z_s)=0$ для всех
     $i$ и
     $s$. Так как
     $\alpha_i\neq 0$, имеем
     $\gamma_s(y_i,z_s)=0$. Таким образом, для любого
     $s$, такого что
     $\gamma_s\neq 0$, выполнены равенства
     $(y_i,z_s)=0$ для всех
     $i\in \{1,2,\dots, k\}$, то есть в
     $h'$ с ненулевыми коэффициентами входят только порождающие
     $z_s$, каждая из которых смежна со всеми
     $y_i$ в графе
     $G$.

     Осталось доказать, что перечисленные условия на линейные
     комбинации порождающих из
     $C(g)$ являются также и достаточными. Пусть
     $g=\sum_{i=1}^k \alpha_i y_i$, где
     $y_i$~--- все порождающие, входящие в линейную комбинацию
     $g$ с ненулевыми коэффициентами,
     $H$~--- подграф графа
     $G$, порожденный множеством вершин
     $\{y_1,\dots, y_k\}$ и
     $H_1,H_2,\dots, H_r$~--- компоненты связности графа
     $\overline{H}$. Имеет место представление
     $g=\sum_{s=1}^r g_s$, где для каждого
     $s=1,2,\dots, r$ элемент
     $g_s$~--- это линейная комбинация порождающих, являющихся
     вершинами графа
     $H_s$.

     Пусть теперь
     $$h=\sum_{i=1}^k \beta_{i} y_i+ \sum_{j=1}^{n-k}\gamma_{j} z_j=\sum_{s=1}^r h_s+\sum_{j=1}^{n-k}\gamma_{j} z_j,$$
     причем для любого
     $s$ существуют
     $\lambda_s,\mu_s\in R$ такие что
     $\lambda_sg_s=\mu_sh_s$ и вершины
     $z_j$ смежны со всеми вершинами
     $y_{i}$, если
     $\gamma_j\neq 0$. Имеем
     \begin{equation}\label{lcprod}
       \begin{aligned}[]
          (g,h) =& \left( \sum_{i=1}^k \alpha_{i}y_{i}, \sum_{l=1}^k \beta_{l}y_{l}+\sum_{j=1}^{n-k}\gamma_{j} z_j\right)=\\
                =&  \sum_{1\leqslant l< i\leqslant k} (\alpha_{i}\beta_l-\alpha_{l}\beta_i) (y_{i},y_l)+
                        \sum_{s=1}^k \sum_{j=1}^{n-k}\alpha_s\gamma_{j} (y_s,z_j)
        \end{aligned}
     \end{equation}
     Если
     $y_i$ и
     $y_l$ лежат в одной компоненте связности
     $H_p$, то из того, что
     $g_p$ и
     $h_p$ пропорциональны, следует, что в разложениях этих элементов в линейные комбинации порождающих
     упорядоченные наборы коэффициентов при порождающих, входящих в
     $H_p$, пропорциональны, а значит по лемме
\ref{proporclemma}
     $\alpha_i\beta_l=\alpha_l\beta_i$, следовательно соответствующее слагаемое равно
     $0$. Если
     $y_i$ и
     $y_l$ лежат в разных компонентах связности графа
     $\overline{H}$, то, в частности,
     $y_i$ и
     $y_l$ не смежны в графе
     $\overline{H}$, а значит смежны в графе
     $H$ и, соответственно, в графе
     $G$. Получаем
     $(y_i,y_l)=0$. Наконец, если
     $\gamma_j\neq 0$, то
     $z_j$ смежна со всеми вершинами
     $y_i$ в графе
     $G$. Следовательно,
     $(y_i,z_j)=0$. Получаем, что выражение
(\ref{lcprod}) равно
     $0$, а значит
     $h\in C(g)$. Таким образом, справедлива следующая теорема.
     \begin{ttt}\label{lincombgen}
       Пусть
       $g$~--- линейная комбинация порождающих алгебры
       $\mathcal{L}(X;G)$,
       $H$~--- подграф графа
       $G$, порожденный вершинами, коэффициенты при которых в
       $g$ отличны от
       $0$ и пусть
       $H_1,H_2,\dots, H_r$~--- компоненты связности графа
       $\overline{H}$. Пусть также
       $g=\sum_{i=1}^r g_i$, где
       $g_i$~--- линейная комбинация порождающих, входящих в граф
       $H_i$. Тогда, в
       $C(g)$ входят линейные комбинации порождающих вида
       $h=\sum_{i=1}^r h_i +h'$, такие что
       $g_i\backsim h_i$ для каждого
       $i\in \{1,2,\dots, r\}$ и
       $h'$~--- линейная комбинация с ненулевыми коэффициентами порождающих, не входящих в
       $H$ и смежных в графе
       $G$ со всеми порождающими, входящими в подграф
       $H$.
     \end{ttt}

   \section{Общий случай} \label{4}
     Этот раздел мы начнем со следующего утверждения о
     PCLS-словах.

     \begin{llll}\label{nonzero}
       Пусть
       $O$~--- некоторый порядок на множестве порождающих
       $X$ алгебры
       $\mathcal{L}(X;G)$. Если
       $[u]$ и
       $[v]$~--- два различных
       $O$-PCLS-слова, начинающиеся с одной и той же буквы, то
       $([u],[v])\neq 0$. Более точно, если
       $[u]>[v]$, то представление
       $([u],[v])$ в виде линейной комбинации
       $O$-базисных мономов имеет вид
       $([u],[v])=[uv]+\sum_{i}\alpha_i [w_i]$, где
       $[w_i]<[uv]$ для всех мономов
       $[w_i]$.
     \end{llll}
     \begin{proof}
       Как следует из леммы~%
\ref{LSAL},
       $[uv]$ является
       $O$-LS-словом. Более того, по лемме~%
\ref{less}, в алгебре
       $\mathcal{L}(X)$ имеет место представление
       \begin{equation} \label{proddecomp}
         ([u],[v])=[uv]+\sum_{i}\alpha_i [w_i],
       \end{equation}
       где
       $[w_i]$~--- различные
       $O$-базисные лиевские мономы алгебры
       $\mathcal{L}(X)$, такие что
       $[w_i]<[uv]$.

       Докажем, что
       $[uv]$ является
       $O$-PCLS-словом. Предположим противное, тогда из
\cite{por11} следует, что ассоциативное слово
       $uv$ содержит подслово
       $w=\tilde{u}y$, являющееся
       $O$-LSA подсловом и такое что
       $y$~--- вторая по старшинству буква слова
       $w$ и
       $y\leftrightarrow \supp ([\tilde{u}])$. При этом старшая буква
       $w$ имеет единственное вхождение в это слово.

       Так как
       $[u]$ и
       $[v]$~---
       $O$-PCLS-слова,
       $w$ не может являться ни подсловом
       $u$, ни подсловом
       $v$. Однако, ситуация, когда слово
       $w$ содержит конец слова
       $u$ и начало слова
       $v$ также невозможна. Действительно, в противном случае первая буква
       $v$ не является первой буквой
       $w$, а так как
       $v$ начинается со старшей буквы слова
       $uv$, слово
       $w$ содержит более одного вхождения этой буквы (так как
       $w$ начинается со старшей своей буквы), получили
       противоречие. Следовательно,
       $[uv]$ является
       $O$-PCLS-словом и, в частности,
       $[uv]\neq 0$ в алгебре
       $\mathcal{L}(X;G)$.

       Так как каждый моном
       $[w_i]$ в алгебре
       $\mathcal{L}(X;G)$ представляется в виде линейной комбинации
       $O$-PCLS-слов, не больших
       $[w_i]$, преобразовав те из мономов
       $[w_i]$ в разложении
(\ref{proddecomp}), которые не являются
       $O$-PCLS-словами, получаем представление
       $([u],[v])=[uv]+h$ в алгебре
       $\mathcal{L}(X;G)$, где
       $h$~--- линейная комбинация
       $O$-PCLS слов, меньших монома
       $[uv]$. Следовательно,
       $[uv]$ не может сократиться с мономами, входящими в
       $h$, откуда и следует утверждение леммы.
     \end{proof}

     \begin{llll}\label{homogennotdisjointsupp}
       Пусть
       $g$ и
       $h$~--- два ненулевых однородных элемента, такие что
       $\supp(g)\cap\, \supp(h)\neq \varnothing$. Тогда
       $(g,h)=0$ тогда и только тогда, когда
       $g\backsim h$.
     \end{llll}
     \begin{proof}
       Пусть
       $x\in \supp(g) \cap\, \supp(h)$. Введем на множестве
       $X$ порядок
       $O$, такой что
       $x$~--- наибольшая буква множества
       $X$.

       Пусть
       $g=\sum_{i=1}^m \alpha_i [u_i]$ и
       $h=\sum_{j=1}^p \beta_j [v_j]$. Без ограничения общности можно считать, что
       все
       $\alpha_i$ и
       $\beta_j$ отличны от нуля,
       $[u_1]<[u_2]<\dots < [u_m]$ и
       $[v_1]<[v_2]<\dots < [v_p]$. Возможны два случая.

       \noindent
       1. Пусть
       $[u_m]\neq [v_p]$. Без ограничения общности можно считать, что
       $[u_m]>[v_p]$. Имеем
       $$(g,h)=\left( \sum_{i=1}^m \alpha_i [u_i], \sum_{j=1}^p \beta_j [v_j]\right)=
           \sum_{i=1}^m \sum_{j=1}^p \alpha_i\beta_j
           ([u_i],[v_j])$$

       Применяя тождество антикоммутативности в тех слагаемых
       $([u_i],[v_j])$, где
       $[v_j]>[u_i]$, получаем
       \begin{equation}\label{homogenprod1}
         (g,h)=\sum_{i=1}^m \sum_{j=1}^p \gamma_{i,j} ([w_{1,i,j}],[w_{2,i,j}]),
       \end{equation}
       где
       $[w_{1,i,j}]=[u_i]$,
       $[w_{2,i,j}]=[v_j]$ и
       $\gamma_{i,j}=\alpha_i\beta_j$, если
       $[u_i]>[v_j]$ и
       $[w_{1,i,j}]=[v_j]$,
       $[w_{2,i,j}])=[u_i]$; и
       $\gamma_{i,j}=-\alpha_i\beta_j$, если
       $[u_i]<[v_j]$.

       По лемме~%
\ref{nonzero} все мономы
       $[w_{1,i,j}w_{2,i,j}]$ являются
       $O$-PCLS-словами, причем
       \begin{equation} \label{PCLSdecom1}
         ([w_{1,i,j}],[w_{2,i,j}])=[w_{1,i,j}w_{2,i,j}]+f_{i,j},
       \end{equation}
       где
       $h_{i,j}$~--- линейная комбинация
       $O$-базисных мономов алгебры
       $\mathcal{L}(X;G)$, каждый из которых меньше
       $[w_{1,i,j}w_{2,i,j}]$.

       Докажем, что наибольшим среди всех мономов вида
       $[w_{1,i,j}w_{2,i,j}]$ является моном
       $[w_{1,m,p}w_{2,m,p}]=[u_m v_p]$. Действительно,
       если
       $[w_{1,i,j}w_{2,i,j}]=[u_i v_j]$, то либо
       $i<m$, откуда
       $u_i<u_m$, а значит
       $u_iv_j<u_mv_p$; либо
       $i=m$, но тогда
       $j<p$, следовательно
       $v_j<v_p$ и получаем
       $u_mv_j<u_mv_p$. Если же
       $[w_{1,i,j}w_{2,i,j}]=[v_j u_i]$, то из того, что
       $u_mv_p\in LSA(X;O)$, следует
       $u_m v_p>v_p u_m>v_ju_i$ (второе неравенство доказывается, как в предыдущем
       случае).

       Таким образом, подставляя
(\ref{PCLSdecom1}) в
 (\ref{homogenprod1}), получаем
       $$(g,h)=\sum_{i=1}^m\sum_{j=1}^p \gamma_{i,j}([w_{1,i,j}w_{2,i,j}]+f_{i,j})=
         \sum_{i=1}^m\sum_{j=1}^p \gamma_{i,j}[w_{1,i,j}w_{2,i,j}]+ \tilde{f},
       $$
       где
       $\tilde{f}=\sum_{i=1}^m\sum_{j=1}^p \gamma_{i,j} f_{i,j}$, а
       значит, по доказанному выше, каждый моном, входящий в
       $\tilde{f}$ заведомо меньше
       $[u_mv_p]$. В итоге получаем
       $(g,h)=\alpha_m\beta_p [u_mv_p]+f$, где
       $f$~--- линейная комбинация
       $O$-PCLS-слов, меньших, чем
       $[u_mv_p]$. Так как
       $\alpha_m\beta_p\neq 0$, следовательно
       $(g,h)\neq 0$.\\

       \noindent
       2. Пусть теперь
       $[u_m]=[v_p]$. Рассмотрим элемент
       $\alpha_m h-\beta_p g$. Если он равен
       $0$ в алгебре
       $\mathcal{L}(X;G)$, то
       $g\backsim h$.

       Пусть
       $\alpha_m h-\beta_p g\neq 0$. Докажем, что в этом случае
       $(g,h)\neq 0$. Пусть это неверно, тогда имеем
       $(g,\alpha_m h-\beta_p g)=\alpha_m (g,h)-\beta_p (g,g)=0$.
       С другой стороны,
       \begin{equation*}
         \begin{aligned}
           \alpha_m h-\beta_p g & =\alpha_m \biggl(\sum_{j=1}^p \beta_j[v_j]\biggr)-\beta_p \biggl(\sum_{i=1}^m \alpha_i[u_i]\biggr)=\\
                                & =\alpha_m \biggl(\sum_{j=1}^{p-1} \beta_j[v_j]\biggr)+\alpha_m \beta_p [v_p]
                                  -\beta_p \biggl(\sum_{i=1}^{m-1}\alpha_i[u_i]\biggr)-\alpha_m \beta_p[u_m]=\\
                                & =\alpha_m \biggl(\sum_{j=1}^{p-1} \beta_j[v_j]\biggr)-\beta_p \biggl(\sum_{i=1}^{m-1}\alpha_i[u_i]\biggr).
         \end{aligned}
       \end{equation*}
       Таким образом, в разложении
       $\alpha_m h-\beta_p g$ в сумму
       $O$-базисных мономов все мономы меньше, чем
       $[u_m]$. По доказанному выше получаем
       $(g,\alpha_m h-\beta_p g)\neq 0$, получили противоречие.
       Следовательно, если
       $g$ и
       $h$ удовлетворяют условиям леммы и не пропорциональны, то
       $(g,h)\neq 0$.

       Докажем обратное утверждение. Пусть
       $g\backsim h$, тогда
       $\alpha g=\beta h$ для некоторых элементов
       $\alpha, \beta \in R$, хотя бы один из которых отличен от
       $0$. Если
       $(g,h)\neq 0$, то, в частности
       $g$ и
       $h$ отличны от
       $0$, а значит
       $\alpha,\beta \in R\backslash \{0\}$. Тогда
       $(\alpha g,\beta h)=\alpha \beta (g,h)\neq 0$, получили
       противоречие с тождеством антикоммутативности в алгебрах
       Ли. Следовательно
       $(g,h)=0$.
     \end{proof}

     \begin{llll}\label{homogendisjointsupp}
       Пусть
       $g$ и
       $h$~--- два ненулевых однородных элемента, такие что
       $\supp(g)\cap\, \supp(h)=\varnothing$. Тогда
       $(g,h)=0$ тогда и только тогда, когда
       $\supp(g)\leftrightarrow \supp(h)$.
     \end{llll}

     \begin{proof}
       Пусть
       $(g,h)=0$ для ненулевых однородных элементов
       $g$ и
       $h$. Фиксируем произвольные буквы
       $x\in \supp (g)$ и
       $y \in \supp (h)$ и определим на множестве
       $X$ порядок
       $O$, в котором
       $x$~--- наибольшая буква, а
       $y$~--- вторая по старшинству.

       Представляя
       $g$ и
       $h$ в виде линейной комбинации
       $O$-PCLS-слов, получаем
       $g=\sum_{i=1}^m \alpha_i [u_i]$ и
       $h=\sum_{j=1}^p \beta_j [v_j]$, где все
       $\alpha_i$ и
       $\beta_j$ отличны от нуля. Отметим, что все мономы
       $[u_i]$ начинаются с
       $x$, а все мономы
       $[v_j]$~--- с
       $y$. Пусть
       $[u_m]$ и
       $[v_p]$~--- наибольшие мономы в разложениях соответственно
       $g$ и
       $h$. Имеем
       \begin{equation} \label{homogenprod3}
         (g,h)=\left(\sum_{i=1}^m \alpha_i [u_i],\sum_{j=1}^p \beta_j[v_j]\right)=
         \sum_{i=1}^m \sum_{j=1}^p \alpha_i\beta_j ([u_i],[v_j]).
       \end{equation}
       Так как
       $[u_i]>_O[v_j]$ для всех
       $i$ и
       $j$, все мономы
       $[u_iv_j]$ являются
       $O$-LS-словами по лемме
\ref{LSAL}, а из леммы
 \ref{less} получаем равенства
       \begin{equation}\label{LSdecom}
         ([u_i],[v_j])=[u_iv_j]+f_{i,j},
       \end{equation}
       где
       $f_{i,j}$~--- линейная комбинация
       $O$-LS-слов меньших, чем
       $[u_iv_j]$. Очевидно, что наибольшим среди мономов
       $[u_iv_j]$ является
       $[u_mv_p]$. Подставляя равенства
(\ref{LSdecom}) в
 (\ref{homogenprod3}), получаем разложение
       $(g,h)$ в линейную комбинацию
       $O$-LS-слов.
       \begin{equation*}
         (g,h)= \sum_{i=1}^m \sum_{j=1}^p \alpha_i\beta_j([u_iv_j]+f_{i,j})
              = \alpha_m\beta_p [u_mv_p]+f,
       \end{equation*}
       где
       $f$~--- линейная комбинация
       $O$-LS-слов, заведомо меньших, чем
       $[u_mv_p]$.

       Если
       $[u_mv_p]$~---
       $O$-PCLS-слово, то оно не может сократиться с мономами,
       входящими в
       $f$. Следовательно,
       $[u_mv_p]$ не является
       $O$-PCLS-словом. Из
\cite{por11} следует, что
       $u_mv_p$ содержит подслово вида
       $w=\tilde{w}z$, являющееся
       $O$-LSA-словом, такое что
       $z$~--- вторая по старшинству буква из входящих в
       $w$ и
       $z\leftrightarrow \supp([\tilde{w}])$. Так как
       $[u_m]$ и
       $[v_p]$~---
       $O$-PCLS-слова, слово
       $w$ не является подсловом ни
       $u_m$, ни
       $v_p$. Следовательно,
       $y\in \supp([w])$, причем эта буква не является первой
       буквой в
       $w$. Отсюда следует, что
       $y$~--- последняя буква слова
       $w$ (и вторая по старшинству буква этого слова), а значит
       $x\in \supp([w])$. Получаем
       $x\leftrightarrow y$.

       В силу произвольности выбора букв
       $x\in \supp(g)$ и
       $y\in \supp(h)$, как наибольшей и второй по старшинству букв множества
       $X$, получаем если
       $(g,h)=0$, то
       $\supp(g)\leftrightarrow \supp(h)$.

       Для доказательства обратного утверждения достаточно показать,
       что если
       $\supp([u])\leftrightarrow \supp([v])$ для
       $O$-PCLS-мономов
       $[u]$ и
       $[v]$ (при произвольном упорядочении
       $O$ множества
       $X$), то
       $([u],[v])=0$.

       Если
       $[v]=y$ для некоторого
       $y\in X$, такого что
       $y\leftrightarrow \supp([u])$, то
       $([u],y)=0$ согласно следствию
\ref{centlincomb}.

       Чтобы доказать, что
       $([u],[v])=0$, если
       $\supp([u])\leftrightarrow \supp([v])$ для произвольных
       $O$-PCLS-мономов, применим индукцию по
       $\ell([v])$. База индукции
       (случай
       $\ell([v])=1$) следует из отмеченного выше. Если
       $\ell([v])>1$, то
       $[v]=([v_1],[v_2])$ для некоторых
       $O$-PCLS-слов
       $[v_1]$ и
       $[v_2]$. Имеем
       \begin{equation} \label{02}
         ([u],[v])=([u],([v_1],[v_2])) = (([u],[v_1]),[v_2])+([v_1],([u],[v_2])).
       \end{equation}
       Из того, что
       $\supp([u]) \leftrightarrow \supp([v])$, следует
       $\supp([u])\leftrightarrow \supp([v_1])$ и
       $\supp([u])\leftrightarrow \supp([v_2])$. Таким образом, по предположению индукции
       $([u],[v_1])=([u],[v_2])=0$. Подставляя эти равенства в
(\ref{02}), получаем
       $([u],[v])=0$.
     \end{proof}

     Фиксируем элемент
     $g\in\mathcal{L}(X;G)$. Пусть
     $H=\overline{G}(\supp(g))$. Обозначим через
     $H_1,H_2,\dots, H_{p}$~--- компоненты связности этого графа.
     \begin{llll}\label{connectcompdecomp}
       Для любого элемента
       $g\in \mathcal{L}(X;G)$ имеет место разложение
       $g=\sum_{i=1}^p g_i$, где
       $g_i$~--- элемент алгебры
       $\mathcal{L}(X;G)$, такой что
       $\supp(g_i)$~--- множество вершин графа
       $H_i$.
     \end{llll}
     \begin{proof}
       Фиксируем произвольное упорядочение
       $O$ множества
       $X$. Нам достаточно доказать, что если
       $[u]$~--- моном, такой что
       $\supp([u])\subseteq \supp(g)$, то
       $[u]$ является произведением порождающих, лежащих только в одной
       компоненте связности графа
       $H$.

       Применим индукцию по
       $\ell([u])$. Если
       $[u]$~--- порождающий, то доказывать нечего. Если
       $\ell([u])=2$, то
       $[u]=(x_i,x_j)$, где
       $x_i>_O x_j$. Если бы
       $x_i$ и
       $x_j$ лежали в различных компонентах связности графа
       $H$, то
       $(x_i,x_j)=0$, то есть моном
       $[u]$ не являлся бы
       $O$-PCLS-словом. Следовательно,
       $x_i$ и
       $x_j$ лежат в одной компоненте связности графа
       $H$.

       Пусть утверждение верно для любого
       $O$-PCLS-монома
       из разложения элемента
       $g$ длина которого не превосходит
       $k\geqslant 3$. Рассмотрим моном
       $[u]$, для которого
       $\ell([u])=k+1$. Имеем
       $[u]=([u_1],[u_2])$, причем
       $\ell([u_j])\leqslant k$ для
       $j=1,2$. По предположению индукции, все порождающие, входящие
       в запись
       $[u_j]$, лежат в некоторой компоненте связности
       $H_{i_j}$ для
       $j=1,2$. Если
       $i_1 \neq i_2$, то, так как вершины, входящие в запись
       $[u_1]$, не соединены с вершинами, входящими в запись
       $[u_2]$, имеем
       $\supp([u_1])\leftrightarrow \supp([u_2])$. Следовательно,
       из леммы~%
\ref{homogendisjointsupp} получаем, что
       $[u]=([u_1],[u_2])=0$ и
       $[u]$ не является
       $O$-PCLS-словом.

       Представим
       $g$ в виде линейной комбинации
       $O$-базисных мономов. В силу доказанного выше, каждый моном
       представляет из себя произведение порождающих из какой-то
       одной компоненты связности графа
       $H$. Таким образом, чтобы получить
       $g_i$, необходимо взять те слагаемые из разложения
       $g$ в линейную комбинацию
       $O$-базисных мономов, которые содержат мономы, в запись
       которых входят только порождающие из множества
       $V(H_i)$
       ($i=1,2,\dots, p$).
     \end{proof}
     \begin{llll} \label{part1}
       Пусть
       $g\in \mathcal{L}(X;G)$ и пусть
       $g=\sum_{i=1}^p g_i$, где
       $g_i$ определены как в лемме~%
\ref{connectcompdecomp}. Тогда
       $C(g)\cap\mathcal{L}(\supp(g),G(\supp(g)))$ состоит из элементов
       вида
       $h=\sum_{i=1}^p h_i$, где
       $g_i$ и
       $h_i$ пропорциональны для
       $i=1,2,\dots, p$.
     \end{llll}
     \begin{proof}
       Сначала докажем, что если
       $h=\sum_{i=1}^p h_i$, где
       $g_i$ и
       $h_i$ пропорциональны для
       $i=1,2,\dots, p$, то
       $h\in C(g)$, а значит и
       $h \in C(g)\cap\mathcal{L}(\supp(g),G(\supp(g)))$.

       Определим графы
       $H,H_1,\dots, H_p$, как и выше. Пусть
       $\widetilde{H}_i$~--- полный граф, множество вершин которого
       совпадает с множеством вершин графа
       $H_i$ для
       $i=1,2,\dots ,p$ и пусть
       $\widetilde{H}$~--- это объединение графов
       $\widetilde{H}_1,\dots, \widetilde{H}_p$.

       Рассмотрим граф
       $\overline{\widetilde{H}}$. В нем порождающие
       $x$ и
       $y$ соединены ребром тогда и только тогда, когда они
       являются вершинами различных графов
       $\widetilde{H}_i$ и
       $\widetilde{H}_j$. Отсюда следует, что
       $\mathcal{L}(\supp(g),\overline{\widetilde{H}})=\bigoplus_{i=1}^p\mathcal{L}(V(H_i))$.
       Так как в свободной алгебре Ли
       централизатор любого элемента
       $f$ состоит из элементов, пропорциональных
       $f$ и только из них, в прямой сумме свободных алгебр Ли
       централизатором элемента
       $g=\sum_{i=1}^p g_i$, очевидно является множество элементов
       вида
       $h=\sum_{i=1}^p h_i$, где
       $g_i$ и
       $h_i$ пропорциональны для
       $i=1,2,\dots, p$.

       Легко видеть, что тождественное отображение порождающих
       $\varphi: \supp(g) \to \supp (g)$ продолжается до
       отображения
       $\varphi: \mathcal{L}(\supp(g), \overline{\widetilde{H}}) \to \mathcal{L}(\supp(g), \overline{H}),$
       которое  является гомоморфизмом, так как все тождества и отношения, выполненные в алгебре
       $\mathcal{L}(\supp(g), \overline{\widetilde{H}})$, выполняются и в алгебре
       $\mathcal{L}(\supp(g), \overline{H})$
       (этот гомоморфизм также будем обозначать
       $\varphi$). Таким образом, если
       $(g,h)=0$ в
       $\mathcal{L}(\supp(g), \overline{\widetilde{H}})$, то
       $(g,h)=0$ и в
       $\mathcal{L}(\supp(g), \overline{H})$. То есть
       $h\in C(g)$ в алгебре
       $\mathcal{L}(\supp(g), \overline{H})$.

       Докажем теперь, что других элементов в множестве
       $C(g)\cap\mathcal{L}(\supp(g),\overline{H})$ нет. Пусть
       $h\in C(g)\cap\mathcal{L}(\supp(g),\overline{H})$. Как и в лемме~
\ref{connectcompdecomp}, легко получить, что
       $h=\sum_{i=1}^p h_i$, где
       $\supp(h_i)\subseteq V(H_i)$. Нам осталось доказать, что если
       $g_i$ и
       $h_i$ не являются пропорциональными хотя бы для одного
       значения
       $i$, то
       $(g,h)\neq 0$.

       Для начала рассмотрим случай, когда
       $H$~--- связный граф. Имеют место разложения
       $g=\sum_{\overline{\delta}}g_{\overline{\delta}}$ и
       $h=\sum_{\overline{\delta}'}h_{\overline{\delta}'}$ элементов
       $g$ и
       $h$ в суммы ненулевых однородных компонент, где
       $\mdeg(g_{\overline{\delta}})=\overline{\delta}$, а
       $\mdeg(h_{\overline{\delta}'})=\overline{\delta}'$.
       Фиксируем порядок
       $\mathcal{O}\in\Omega(x)$ для произвольного
       $x\in \supp(h)$. Пусть
       $\overline{\varepsilon}$ и
       $\overline{\varepsilon}'$~--- наибольшие среди мультистепеней однородных
       элементов относительно выбранного порядка.

       Отметим, что в этом случае
       $x\in \supp(g_{\overline{\varepsilon}})\cap \supp(h_{\overline{\varepsilon}}')$. Действительно,
       $x\in \supp(g_{\overline{\delta}})$ и
       $x\in \supp(h_{\overline{\delta}'})$ хотя бы для одной
       мультистепени
       $\overline{\delta}$ и хотя бы для одной мультистепени
       $\overline{\delta}'$; в силу выбора порядка на
       мультистепенях получаем требуемое.

       Возможны два случая.\\

       \noindent
       1. Пусть
       $g_{\overline{\varepsilon}}$ и
       $h_{\overline{\varepsilon}'}$ не пропорциональны.  Тогда из леммы~
\ref{homogennotdisjointsupp} следует, что
       $(g_{\overline{\varepsilon}},h_{\overline{\varepsilon}'})\neq 0$.

       Имеем
       $(g,h)=\sum_{\overline{\delta}:\overline{\delta}'}(g_{\overline{\delta}},h_{\overline{\delta}'})$.
       Из только что доказанного следует, что в разложении
       $(g,h)$ в сумму однородных элементов есть ненулевое
       слагаемое
       $(g_{\overline{\varepsilon}},h_{\overline{\varepsilon}'})$ мультистепени
       $\overline{\varepsilon}+\overline{\varepsilon}'$, причем мультистепени всех остальных слагаемых меньше
       относительно порядка
       $\mathcal{O}$. Следовательно,
       $(g_{\overline{\varepsilon}},h_{\overline{\varepsilon}'})$
       не может сократиться с другими слагаемыми, входящим в
       разложение
       $(g,h)$, а значит
       $(g,h)\neq 0$.\\

       \noindent
       2. Теперь предположим, что
       $g_{\overline{\varepsilon}}$ и
       $h_{\overline{\varepsilon}'}$ пропорциональны, то есть
       $\overline{\varepsilon}=\overline{\varepsilon}'$ и
       $\alpha g_{\overline{\varepsilon}}= \beta h_{\overline{\varepsilon}}$ для некоторых
       $\alpha,\beta \in R$, хотя бы одно из которых не равно
       $0$ (отметим, что так как
       $g_{\overline{\varepsilon}}$ и
       $h_{\overline{\varepsilon}'}$ не равны нулю, получаем
       $\alpha, \beta \in R\backslash \{0\}$).

       Пусть
       $(g,h)=0$. Рассмотрим элемент
       $f=\beta h-\alpha g$. Если он равен
       $0$, то
       $g$ и
       $h$ пропорциональны. Предположим, что
       $\beta h-\alpha g\neq 0$. Тогда имеем
       \begin{equation}\label{prodzero}
         (g,\beta h-\alpha g)=\beta(g,h)-\alpha(g,g)=0.
       \end{equation}
       С другой стороны,
       \begin{equation} \label{perform}
         \begin{aligned}
           \beta h- \alpha g & =\beta \sum_{\overline{\delta}'} h_{\overline{\delta}'}-
           \alpha \sum_{\overline{\delta}}g_{\overline{\delta}}\\
           & =\beta\biggl(h_{\overline{\varepsilon}}+ \sum_{\overline{\delta}':\overline{\delta}'< \overline{\varepsilon}} h_{\overline{\delta}'}\biggr)
            - \alpha\biggl(g_{\overline{\varepsilon}}+ \sum_{\overline{\delta}:\overline{\delta}< \overline{\varepsilon}}g_{\overline{\delta}}\biggr)\\
         \end{aligned}
       \end{equation}
       \begin{equation*}
         \begin{aligned}
           & =(\beta h_{\overline{\varepsilon}}- \alpha g_{\overline{\varepsilon}})
              +\beta \sum_{\overline{\delta}':\overline{\delta}'< \overline{\varepsilon}} h_{\overline{\delta}'}
              - \alpha \sum_{\overline{\delta}:\overline{\delta}< \overline{\varepsilon}}g_{\overline{\delta}}\\
           & =\beta \sum_{\overline{\delta}':\overline{\delta}'< \overline{\varepsilon}} h_{\overline{\delta}'}
              - \alpha \sum_{\overline{\delta}:\overline{\delta}< \overline{\varepsilon}}g_{\overline{\delta}}
         \end{aligned}
       \end{equation*}
       Выражение, полученное после приведения подобных в
(\ref{perform}), обозначим
       $\sum_{\overline{\xi}} f_{\overline{\xi}}$ то есть
       $f=\sum_{\overline{\xi}} f_{\overline{\xi}}$. Таким образом,
       мультистепени всех слагаемых в разложении
       $f$ меньше
       $\overline{\varepsilon}$ относительно порядка
       $\mathcal{O}$.  Пусть
       $f_{\overline{\zeta}}$~--- моном наибольшей мультистепени из разложения
       $f$ в сумму однородных мономов.  Разберем возможные случаи.\\

       \noindent
       2-1. Если
       $x\in \supp(f)$, то
       $x\in \supp (g_{\overline{\epsilon}})\cap \supp (f_{\overline{\zeta}})$,
       и мы оказались в ситуации, рассмотренной в случае 1. Таким образом,
       $(g,f)\neq 0$, получили противоречие с
(\ref{prodzero}).\\

       \noindent
       Пусть теперь
       $x\not \in \supp(f)$ и пусть
       $(g_{\overline{\delta}},f_{\overline{\xi}})\neq 0$ для
       некоторых
       $g_{\overline{\delta}}$ и
       $f_{\overline{\xi}}$. Обозначим через
       $\overline{\delta}_1$ наибольшую мультистепень (относительно порядка
       $\mathcal{O}$), такую что
       $(g_{\overline{\delta}_1},f_{\overline{\xi}})\neq 0$ для
       некоторой мультистепени
       $\overline{\xi}$ . Наибольшую такую мультистепень в разложении
       элемента
       $f$ обозначим
       $\overline{\xi}_1$.\\

       \noindent
       2-2. Пусть
       $x\in \supp (g_{\overline{\delta}_1})$. Данный случай снова разбивается на
       два.\\

       \noindent
       2-2-1. Предположим, что
       $\supp (g_{\overline{\delta}_1})\cap \supp(f_{\overline{\xi}_1}) \neq \varnothing$
       и пусть
       $y$~--- произвольный элемент этого множества.

       Введем порядок
       $\mathcal{O}_1\in\Omega(x,y)$ (порядки
       $\mathcal{O}$ и
       $\mathcal{O}_1$ могут быть различными).

       Пусть
       $\overline{\delta}_0$~--- наибольшая (в смысле порядка
       $\mathcal{O}_1$) из мультистепеней мономов, входящих в разложение
       $g$, не больших мультистепени
       $\overline{\delta}_1$ (в смысле порядка
       $\mathcal{O}$), а
       $\overline{\xi}_0$~--- наибольшая из мультистепеней
       мономов, входящих в разложение
       $f$. Тогда
       $\overline{\delta}_0 \geqslant_{\mathcal{O}_1} \overline{\delta}_1$.
       Таким образом, из
       $x,y\in \supp (g_{\overline{\delta}_1})$ следует
       $x,y\in \supp (g_{\overline{\delta}_0})$. Аналогично,
       $y \in \supp (f_{\overline{\xi}_0})$. Так как
       $x\not \in \supp (f_{\overline{\xi}_0})$, имеем
       $\overline{\delta}_0\neq \overline{\xi}_0$. Следовательно,
       $(g_{\overline{\delta}_0},f_{\overline{\xi}_0})\neq 0$.
       Кроме того, все остальные однородные слагаемые в
       разложении произведения
       $(g,f)$ имеют меньшую мультистепень относительно порядка
       $\mathcal{O}_1$. Действительно,
       $f_{\overline{\xi}_0}$ имеет наибольшую степень, среди
       слагаемых, входящих в разложение
       $f$, а все слагаемые
       $g_{\overline{\delta}}$ при
       $\overline{\delta}>\overline{\delta}_0$ удовлетворяют
       условию
       $(g_{\overline{\delta}},f_{\overline{\xi}})=0$ для любого
       $f_{\overline{\xi}}$. Таким образом, элемент
       $(g_{\overline{\delta}_0},f_{\overline{\xi}_0})$ не может
       сократиться с другими элементами, входящими в разложение
       $(g,f)$, так как их мультистепени меньше, чем
       $\overline{\delta}_0+\overline{\xi}_0$ в смысле порядка
       $\mathcal{O}_1$. Следовательно,
       $(g,f)\neq 0$, получили противоречие с
(\ref{prodzero}).\\

       \noindent
       2-2-2.
       Предположим, что
       $\supp (g_{\overline{\delta}_1})\cap \supp(f_{\overline{\xi}_1})=\varnothing$.
       Тогда существуют порождающие
       $y \in \supp(g_{\overline{\delta}_1})$ и
       $z \in \supp(f_{\overline{\xi}_1})$, смежные в графе
       $H$.

       Рассмотрим на множестве мультистепеней порядок
       $\mathcal{O}_1\in \Omega(x,y,z)$. Пусть
       $\overline{\delta}_0$~--- наибольшая (в смысле порядка
       $\mathcal{O}_1$) из мультистепеней мономов,
       входящих в разложение
       $g$, не больших мультистепени
       $\overline{\delta}_1$ (в смысле порядка
       $\mathcal{O}$), а
       $\overline{\xi}_0$~--- наибольшая из мультистепеней
       мономов, входящих в разложение
       $f$ (опять относительно порядка
       $\mathcal{O}_1$). Имеем
       $\overline{\delta}_0 \geqslant_{\mathcal{O}_1} \overline{\delta}_1$,
       поэтому из
       $x,y\in \supp (g_{\overline{\delta}_1})$ следует
       $x,y\in \supp (g_{\overline{\delta}_0})$. Так как
       $z \in \supp (f_{\overline{\xi}_1})$, получаем либо
       $z \in \supp (f_{\overline{\xi}_0})$, либо
       $y \in \supp (f_{\overline{\xi}_0})$. Так как
       $x\not \in \supp (f_{\overline{\xi}_0})$, имеем
       $\overline{\delta}_0\neq \overline{\xi}_0$. Следовательно,
       в обоих случаях
       $(g_{\overline{\delta}_0},f_{\overline{\xi}_0})\neq 0$ (в первом случае, в силу
       леммы~%
\ref{homogendisjointsupp}, во втором~--- в силу леммы
 \ref{homogennotdisjointsupp}).
       Кроме того, как и в случае 2-2-1,
       все остальные однородные слагаемые в
       разложении произведения
       $(g,f)$ имеют меньшую мультистепень относительно порядка
       $\mathcal{O}_1$. Таким образом, элемент
       $(g_{\overline{\delta}_0},f_{\overline{\xi}_0})$ не может
       сократиться с другими элементами, входящими в разложение
       $(g,f)$. Следовательно,
       $(g,f)\neq 0$, и мы снова получили противоречие с
(\ref{prodzero}).\\

       \noindent
       2-3. Теперь рассмотрим случай
       $x \not \in \supp (g_{\overline{\delta}_1})$. Имеет место представление
       $g=g'+g''$, где
       $g'$~--- сумма всех слагаемых
       $g_{\overline{\delta}}$, таких что
       $\supp(g_{\overline{\delta}})\leftrightarrow \supp(f)$, а
       $g''$~--- сумма всех остальных слагаемых. Имеем
       $(g_{\overline{\delta}},f_{\overline{\xi}})=0$ для любого элемента
       $g_{\overline{\delta}})$, входящего в разложение
       $g'$, и для любого элемента
       $f_{\overline{\xi}}$, входящего в разложение
       $f$. Кроме
       того, так как
       $\supp(f) \subseteq \supp(g)$ имеем
       $g''\neq 0$. Получаем два случая.\\

       \noindent
       2-3-1. Предположим, что
       $\supp(g') \nleftrightarrow \supp(g'')$. Если
       $\supp(g')\cap \supp(g'')\neq \varnothing$, то через
       $y$ обозначим произвольный элемент этого множества. Если же
       $\supp(g')\cap \supp(g'')=\varnothing$, то пусть
       $y \in \supp(g'')$~--- это элемент, такой что
       $y \nleftrightarrow \supp(g')$. Пусть
       $\mathcal{O}_1=\Omega(x,y)$ и пусть
       $g_{\overline{\delta}_2}$~--- наибольший элемент из разложения
       $g''$ относительно этого порядка, такой что
       $y\in \supp(g_{\overline{\delta}_2})$. Снова нужно рассмотреть
       два случая.\\

       \noindent
       2-3-1-1. Пусть
       $\supp(g_{\overline{\delta}_2})\cap \supp (f) \neq \varnothing$, то есть существует элемент
       $f_{\overline{\xi}_1}$ из разложения
       $f$ в сумму однородных элементов, такой что
       $\supp(g_{\overline{\delta}_2})\cap \supp (f_{\overline{\xi}_1}) \neq \varnothing$ и
       пусть
       $z$~--- произвольный порождающий из этого множества.

       Пусть
       $\mathcal{O}_2\in \Omega(x,y,z)$ и пусть
       $\overline{\delta}_0$~--- наибольшая (в смысле этого порядка) из мультистепеней мономов,
       входящих в разложение
       $g''$, а
       $\overline{\xi}_0$~--- наибольшая из мультистепеней
       мономов, входящих в разложение
       $f$ (опять в смысле порядка
       $\mathcal{O}_2$). Имеем
       $\overline{\delta}_0 \geqslant_{\mathcal{O}_2} \overline{\delta}_2$, а так как
       $\overline{\delta}_0 \leqslant_{\mathcal{O}_1} \overline{\delta}_2$, степени относительно
       $y$ у
       $\overline{\delta}_0$ и
       $\overline{\delta}_2$ одинаковы и степень относительно
       $z$ у
       $\overline{\delta}_0$ не меньше степени относительно
       $z$ у
       $\overline{\delta}_2$. Значит из
       $y,z\in \supp (g_{\overline{\delta}_2})$ следует
       $y,z\in \supp (g_{\overline{\delta}_0})$. Аналогично,
       $z \in \supp (f_{\overline{\xi}_0})$. А так как
       $y\nleftrightarrow \supp(g')$ и
       $\supp(f)\leftrightarrow \supp(g')$, получаем
       $y \not \in \supp (f_{\overline{\xi}_0})$, следовательно,
       $\overline{\delta}_0\neq \overline{\xi}_0$. Таким образом,
       $(g_{\overline{\delta}_0},f_{\overline{\xi}_0})\neq 0$ по лемме~%
\ref{homogennotdisjointsupp}.
       Кроме того, аналогично предыдущим случаям доказывается, что все остальные однородные слагаемые в
       разложении произведения
       $(g,f)$ имеют меньшую мультистепень относительно нового порядка. Таким образом, элемент
       $(g_{\overline{\delta}_0},f_{\overline{\xi}_0})$ не может
       сократиться с другими элементами, входящими в разложение
       $(g,f)$. Следовательно,
       $(g,f)\neq 0$, получили противоречие с
(\ref{prodzero}).\\

       \noindent
       2-3-1-2. Пусть
       $\supp(g_{\overline{\delta}_2})\cap \supp (f)=\varnothing$, но
       $\supp(g_{\overline{\delta}_2})\nleftrightarrow \supp (f)$, то есть существует элемент
       $f_{\overline{\xi}_1}$ из разложения
       $f$ в сумму однородных элементов, такой что для некоторых
       порождающих
       $z\in \supp(g_{\overline{\delta}_2})$ и
       $t\in \supp (f_{\overline{\xi}_1}) \neq 0$  выполнено
       $z \nleftrightarrow t$.

       Пусть
       $\mathcal{O}_2\in \Omega(x,y,z,t)$ и пусть
       $\overline{\delta}_0$~--- наибольшая (в смысле этого порядка) из мультистепеней мономов,
       входящих в разложение
       $g''$, а
       $\overline{\xi}_0$~--- наибольшая (опять в смысле порядка
       $\mathcal{O}_2$) из мультистепеней
       мономов, входящих в разложение
       $f$. Имеем
       $\overline{\delta}_0 \geqslant_{\mathcal{O}_2} \overline{\delta}_2$, а так как
       $\overline{\delta}_0 \leqslant_{\mathcal{O}_1} \overline{\delta}_2$, степени относительно
       $y$ у
       $\overline{\delta}_0$ и
       $\overline{\delta}_2$ одинаковы и степень относительно
       $z$ у
       $\overline{\delta}_0$ не меньше степени относительно
       $z$ у
       $\overline{\delta}_2$.
       Таким образом, из
       $y,z\in \supp (g_{\overline{\delta}_2})$ следует
       $y,z\in \supp (g_{\overline{\delta}_0})$.  Так как
       $y,z \not \in \supp (f_{\overline{\xi}_0})$, степень
       относительно
       $t$ у
       $\overline{\xi}_0$ не меньше степени относительно
       $t$ у
       $\overline{\xi}_1$, следовательно,
       $t\in \supp(f_{\overline{\xi}_0})$. В итоге
       получаем, что
       $\overline{\delta}_0\neq \overline{\xi}_0$. Следовательно,
       $(g_{\overline{\delta}_0},f_{\overline{\xi}_0})\neq 0$ по лемме~%
\ref{homogennotdisjointsupp}.
       Кроме того, аналогично предыдущим случаям доказывается, что все остальные однородные слагаемые в
       разложении произведения
       $(g,f)$ имеют меньшую мультистепень относительно порядка
       $\mathcal{O}_1$. Таким образом, элемент
       $(g_{\overline{\delta}_0},f_{\overline{\xi}_0})$ не может
       сократиться с другими элементами, входящими в разложение
       $(g,f)$. Следовательно,
       $(g,f)\neq 0$, получили противоречие с
(\ref{prodzero}).\\

       \noindent
       2-3-2. Пусть
       $\supp(g') \leftrightarrow \supp(g'')$. То есть в графе
       $H$ вершины из множества
       $\supp(g')$ не смежны с вершинами из множества
       $\supp(g'')$. Получили противоречие со связностью графа
       $H$.\\

       \noindent
       2-4. Пусть все произведения
       $(g_{\overline{\delta}},f_{\overline{\xi}})$ равны
       $0$. Рассмотрим произвольный элемент
       $f_{\overline{\xi}_0}$ из разложения
       $f$ в сумму однородных элементов. Так как
       $\supp(f)\subseteq \supp (g)$, имеем
       $\supp(f_{\overline{\xi}_0}) \subseteq \supp (g)$. То есть
       $\supp (f_{\overline{\xi}_0})\cap \supp(g_{\overline{\delta}_0})\neq \varnothing$
       для некоторого слагаемого
       $g_{\overline{\delta}_0}$. Так как
       $(g_{\overline{\delta}_0},f_{\overline{\xi}_0})=0$, из
       леммы~%
\ref{homogennotdisjointsupp} следует, что
       $g_{\overline{\delta}_0}$ и
       $f_{\overline{\xi}_0}$ пропорциональны. Следовательно,
       $(g_{\overline{\delta}_0}, g_{\overline{\delta}})=0$ для
       всех
       ${\overline{\delta}}$ из разложения
       $g$ в сумму однородных элементов. Из лемм~%
\ref{homogennotdisjointsupp} и
 \ref{homogendisjointsupp} следует, что
       $\supp(g_{\overline{\delta}_0})\leftrightarrow \supp(g-g_{\overline{\delta}_0})$,
       следовательно, вершины множеств
       $\supp(g_{\overline{\delta}_0})$ не связаны с вершинами
       множества
       $\supp(g) \backslash \supp(g_{\overline{\delta}_0})$, что
       противоречит связности графа
       $H$.\\

       Перейдем к рассмотрению общего случая. Рассмотрим
       произвольный элемент
       $g$. И пусть
       $h\in C(g)\cap \mathcal{L}(\supp(g), \overline{H})$.

       Пусть
       $g=\sum_{i=1}^p g_i$ и
       $h=\sum_{j=1}^p h_j$, заданные в условии леммы. Так как
       $\supp(g_i) \leftrightarrow \supp(h_j)$ при
       $i\neq j$, в силу леммы~%
\ref{homogendisjointsupp} имеем
       $0=(g,h)=\sum_{i=1}^p (g_i,h_i)$, а так как
       $\supp(g_i)$ и
       $\supp(h_i)$ состоят из порождающих, являющихся вершинами
       графа
       $H_i$, получаем
       $(g_i,h_i)=0$ для
       $i=1,2,\dots, p$. Из доказанного выше для связных графов
       следует, что
       $g_i$ и
       $h_i$ пропорциональны для
       $i=1,2,\dots, p$, что и требовалось доказать.
     \end{proof}
     \begin{llll}\label{part2}
       Пусть
       $g\in \mathcal{L}(X;G)$ и пусть
       $h\in C(g)$, такой что
       $\supp(\widetilde{h}) \nsubseteq \supp(g)$ ни для какого однородного элемента
       $\widetilde{h}$ из разложения
       $h$. Тогда
       $\supp(h) \leftrightarrow \supp(g)$.
     \end{llll}
     \begin{proof}
       Пусть
       $h=h'+h''$, где
       $h'$ состоит из всех элементов
       $\widetilde{h}$ разложения
       $h$ в сумму однородных элементов, для которых
       $g\leftrightarrow \widetilde{h}$. Имеем
       $(g,h'')=(g,h-h'')=(g,h)-(g,h'')=0$. \\

       \noindent
       1. Пусть
       $\supp(g)\cap \supp(h'')\neq \varnothing$ и пусть
       $x$~--- произвольный элемент из этого множества. Рассмотрим
       порядок
       $\mathcal{O}\in \Omega(x)$.

       Пусть
       $g_{\overline{\varepsilon}}$ и
       $h_{\overline{\eta}}$~--- элементы наибольших мультистепеней относительно
       определенного выше порядка в разложениях
       $g$ и
       $h''$, соответственно. Тогда
       $x \in \supp (g_{\overline{\varepsilon}})\cap \supp(h_{\overline{\eta}})$.
       Так как
       $\supp(h_{\overline{\eta}}) \nsubseteq \supp(g)$, получаем
       $\overline{\varepsilon} \neq \overline{\eta}$, следовательно,
       $g_{\overline{\varepsilon}}$ и
       $h_{\overline{\eta}}$ не являются пропорциональными, а
       значит по лемме~%
\ref{homogennotdisjointsupp}
       $(g_{\overline{\varepsilon}},h_{\overline{\eta}}) \neq 0$.
       Отсюда следует, что наибольшая из степеней слагаемых в разложении
       $(g,h'')$ в сумму однородных элементов равна
       $\overline{\varepsilon}+\overline{\eta}$. Этот элемент не может сократиться
       с другими элементами разложения, так как они имеют меньшую
       мультистепень относительно порядка
       $\mathcal{O}$. Следовательно,
       $(g,h'')\neq 0$.\\

       \noindent
       2. Из доказанного выше следует, что
       $\supp(g)\cap \supp(h'')=\varnothing$. Следовательно,
       существуют
       $x\in \supp(g)$ и
       $y\in \supp(h'')$, такие что
       $x \nleftrightarrow y$.

       Пусть
       $\mathcal{O}\in\Omega(x,y)$ и пусть
       $g_{\overline{\varepsilon}}$ и
       $h_{\overline{\eta}}$~--- элементы наибольших
       мультистепеней в смысле этого порядка, входящие в
       разложения
       $g$ и
       $h''$ соответственно. В силу выбора  этих элементов,
       по лемме~%
\ref{homogendisjointsupp} получаем
       $(g_{\overline{\varepsilon}},h_{\overline{\eta}}) \neq 0$.
       Отсюда следует, что наибольшая из степеней слагаемых в разложении
       $(g,h'')$ в сумму однородных элементов равна
       $\overline{\varepsilon}+\overline{\eta}$. Как и в предыдущем случае, этот элемент не может сократиться
       с другими элементами разложения, так как они имеют меньшую
       мультистепень. Следовательно,
       $(g,h'')\neq 0$, получили противоречие.\\

       Итак, мы получили, что
       $h''=0$. Отсюда следует утверждение леммы.
     \end{proof}
     Теперь мы можем описать множество
     $C(g)$.
     \begin{ttt}\label{gencase}
       Пусть
       $g\in \mathcal{L}(X;G)$,
       $H=\overline{G}(\supp(g))$ и
       $H_1,\dots, H_p$~--- компоненты связности графа
       $H$. Пусть также
       $g=\sum_{i=1}^p g_i$, где
       $\supp(g_i)$ состоит из вершин графа
       $H_i$ для всех
       $i=1,2,\dots p$. Тогда
       $C(g)$ состоит из элеметов вида
       $h=\sum_{i=1}^p h_i+h'$, где
       $g_i\backsim h_i$ для
       $i=1,2,\dots p$ и
       $\supp(g)\leftrightarrow \supp(h')$.
     \end{ttt}
     \begin{proof}
       Пусть
       $h\in C(g)$. Тогда представим
       $h$ в виде
       $h=\check{h}+h'$, где
       $\check{h}$ состоит из всех однородных элементов
       $f$ разложения
       $h$, таких что
       $\supp(f)\subseteq \supp(g)$. Легко видеть, что из
       $(g,h)=0$ следует
       $(g,\check{h})=0$ и
       $(g,h')=0$. Действительно, в противном случае для любого слагаемого
       $f$ из разложения
       $(g,\check{h})$ в сумму однородных элементов имеет
       место включение
       $\supp(f)\subseteq \supp(g)$. С другой стороны, для любого
       слагаемого
       $d$ из разложения
       $(g,h')$ в сумму однородных элементов выполнено
       $\supp(d)\nsubseteq \supp(g)$, а значит слагаемые,
       входящие в
       $(g,f)$, не могут сократиться со слагаемыми из
       $(g,h')$.

       Так как
       $\check{h}\in C(g)\cap \mathcal{L}(\supp(g);\overline{H})$, из
       леммы~%
\ref{part1} получаем
       $\check{h}=\sum_{i=1}^p h_i$, где
       $g_i$ и
       $h_i$ пропорциональны для
       $i=1,2,\dots, p$. А так как
       $h'$ удовлетворяет условию леммы~%
\ref{part2},
       $\supp(g)\leftrightarrow \supp(h')$. Теорема доказана.
     \end{proof}
     
 \end{document}